\documentclass[11pt]{amsart}

\usepackage{amsmath,amssymb,amscd}
\input{prepictex}
\input{pictex}
\input{postpictex}

\begin{document}
\theoremstyle{plain}
\newtheorem{thm}{Theorem}[section]
\newtheorem*{thm*}{Theorem}
\newtheorem{prop}[thm]{Proposition}
\newtheorem*{prop*}{Proposition}
\newtheorem{lemma}[thm]{Lemma}
\newtheorem{cor}[thm]{Corollary}
\newtheorem*{conj*}{Conjecture}
\newtheorem*{cor*}{Corollary}
\newtheorem{defn}[thm]{Definition}
\newtheorem{ocond}{Old Condition}
\newtheorem{ncond}{New Condition}
\theoremstyle{definition}
\newtheorem*{defn*}{Definition}
\newtheorem{rems}[thm]{Remarks}
\newtheorem*{rems*}{Remarks}
\newtheorem*{proof*}{Proof}
\newtheorem*{not*}{Notation}
\newcommand{\npartial}{\slash\!\!\!\partial}
\newcommand{\Heis}{\operatorname{Heis}}
\newcommand{\Solv}{\operatorname{Solv}}
\newcommand{\Spin}{\operatorname{Spin}}
\newcommand{\SO}{\operatorname{SO}}
\newcommand{\ind}{\operatorname{ind}}
\newcommand{\Index}{\operatorname{index}}
\newcommand{\ch}{\operatorname{ch}}
\newcommand{\rank}{\operatorname{rank}}

\newcommand{\abs}[1]{\lvert#1\rvert}
 \newcommand{\A}{{\mathcal A}}
        \newcommand{\D}{{\mathcal D}}\newcommand{\HH}{{\mathcal H}}
        \newcommand{\LL}{{\mathcal L}}
        \newcommand{\B}{{\mathcal B}}
        \newcommand{\K}{{\mathcal K}}
\newcommand{\oo}{{\mathcal O}}
         \newcommand{\PP}{{\mathcal P}}
        \newcommand{\s}{\sigma}
\newcommand{\al}{\alpha}
        \newcommand{\coker}{{\mbox coker}}
        \newcommand{\p}{\partial}
        \newcommand{\dd}{|\D|}
        \newcommand{\n}{\parallel}
\newcommand{\bma}{\left(\begin{array}{cc}}
\newcommand{\ema}{\end{array}\right)}
\newcommand{\bca}{\left(\begin{array}{c}}
\newcommand{\eca}{\end{array}\right)}
\def\clsp{\overline{\operatorname{span}}}
\def\T{\mathbb T}
\def\Aut{\operatorname{Aut}}

\newcommand{\sr}{\stackrel}
\newcommand{\da}{\downarrow}
\newcommand{\tD}{\tilde{\D}}

        \newcommand{\R}{\mathbf R}
        \newcommand{\C}{\mathbf C}
        \newcommand{\h}{\mathbf H}
\newcommand{\Z}{\mathbf Z}
\newcommand{\N}{\mathbf N}
\newcommand{\tto}{\longrightarrow}
\newcommand{\ben}{\begin{displaymath}}
        \newcommand{\een}{\end{displaymath}}
\newcommand{\be}{\begin{equation}}
\newcommand{\ee}{\end{equation}}

        \newcommand{\bean}{\begin{eqnarray*}}
        \newcommand{\eean}{\end{eqnarray*}}
\newcommand{\nno}{\nonumber\\}
\newcommand{\bea}{\begin{eqnarray}}
        \newcommand{\eea}{\end{eqnarray}}

\def\cross#1{\rlap{\hskip#1pt\hbox{$-$}}}
        \def\intcross{\cross{0.3}\int}
        \def\bigintcross{\cross{2.3}\int}

\newcommand{\supp}[1]{\operatorname{#1}}
\newcommand{\norm}[1]{\parallel\, #1\, \parallel}
\newcommand{\ip}[2]{\langle #1,#2\rangle}
\setlength{\parskip}{.3cm}
\newcommand{\nc}{\newcommand}
\nc{\nt}{\newtheorem} \nc{\gf}[2]{\genfrac{}{}{0pt}{}{#1}{#2}}
\nc{\mb}[1]{{\mbox{$ #1 $}}} \nc{\real}{{\mathbb R}}
\nc{\comp}{{\mathbb C}} \nc{\ints}{{\mathbb Z}}
\nc{\Ltoo}{\mb{L^2({\mathbf H})}} \nc{\rtoo}{\mb{{\mathbf R}^2}}
\nc{\slr}{{\mathbf {SL}}(2,\real)} \nc{\slz}{{\mathbf
{SL}}(2,\ints)} \nc{\su}{{\mathbf {SU}}(1,1)} \nc{\so}{{\mathbf
{SO}}} \nc{\hyp}{{\mathbb H}} \nc{\disc}{{\mathbf D}}
\nc{\torus}{{\mathbb T}}
\newcommand{\tk}{\widetilde{K}}
\newcommand{\boe}{{\bf e}}\newcommand{\bt}{{\bf t}}
\newcommand{\vth}{\vartheta}
\newcommand{\CGh}{\widetilde{\CG}}
\newcommand{\db}{\overline{\partial}}
\newcommand{\tE}{\widetilde{E}}
\newcommand{\tr}{\mbox{tr}}
\newcommand{\ta}{\widetilde{\alpha}}
\newcommand{\tb}{\widetilde{\beta}}
\newcommand{\txi}{\widetilde{\xi}}
\newcommand{\hV}{\hat{V}}
\newcommand{\IC}{\mathbf{C}}
\newcommand{\IZ}{\mathbf{Z}}
\newcommand{\IP}{\mathbf{P}}
\newcommand{\IR}{\mathbf{R}}
\newcommand{\IH}{\mathbf{H}}
\newcommand{\IG}{\mathbf{G}}
\newcommand{\CC}{{\mathcal C}}
\newcommand{\CS}{{\mathcal S}}
\newcommand{\CG}{{\mathcal G}}
\newcommand{\CL}{{\mathcal L}}
\newcommand{\CO}{{\mathcal O}}
\nc{\ca}{{\mathcal A}} \nc{\cag}{{{\mathcal A}^\Gamma}}
\nc{\cg}{{\mathcal G}} \nc{\chh}{{\mathcal H}} \nc{\ck}{{\mathcal
B}} \nc{\cl}{{\mathcal L}} \nc{\cm}{{\mathcal M}}
\nc{\cn}{{\mathcal N}} \nc{\cs}{{\mathcal S}} \nc{\cz}{{\mathcal
Z}}
\nc{\sind}{\sigma{\rm -ind}}
\newcommand{\la}{\langle}
\newcommand{\ra}{\rangle}
\newcommand{\cda}{{\mathcal C}_\D(\A)}
\renewcommand{\labelitemi}{{}} 

\newcommand{\field}[1]{\mathbb{#1}}
\newcommand{\NN}{\field{N}}

\newcommand{\Mm}{\mathcal{M}}

\newcommand{\Id}{\operatorname{id}}
\newcommand{\One}{\ensuremath{\mathbf{1}}}
\newcommand{\lsp}{\operatorname{span}}
\begin{center}
 \title{Noncommutative Manifolds from Graph and $k$-Graph $C^*$-Algebras}

 \vspace{.5 in}

\author{David Pask, Adam Rennie, Aidan Sims}\maketitle
\vspace{-0.2in}
{\em School of Mathematical and Physical Sciences\\
University of Newcastle, Callaghan, NSW
Australia, 2308\\ david.pask@newcastle.edu.au,\\
aidan.sims@newcastle.edu.au\\
Institute for Mathematical Sciences, University of Copenhagen\\
Universitetsparken 5, DK-2100, Copenhagen\\ rennie@math.ku.dk}

\end{center}

\centerline{Abstract} In \cite{PRen} we constructed smooth
$(1,\infty)$-summable semfinite spectral triples for graph algebras
with a faithful trace, and in \cite{PRS} we constructed
$(k,\infty)$-summable semifinite spectral triples for $k$-graph
algebras. In this paper we identify classes of graphs and $k$-graphs
which satisfy a version of Connes' conditions for noncommutative
manifolds.

\tableofcontents
\section{Introduction}\label{mflds}

The object of this paper is to address the general definition of
noncommutative manifolds. The phrase `noncommutative manifold' is
one which is still open to some degree of interpretation. Broadly
speaking, a noncommutative manifold is a spectral triple
$(\A,\HH,\D)$ satisfying some additional conditions, such as those
originally proposed by Connes, \cite{C1}. However, the conditions
presented in \cite{C1} only make sense when $\mathcal{A}$ is unital
(that is, a compact noncommutative space). Moreover, the proof of
Connes' spin manifold theorem in \cite{RV} uses a modification of
Connes' conditions even in the compact case.

We consider the set of conditions presented in \cite{RV}, and show
that there is a natural generalisation of each to noncommutative,
nonunital and semifinite spectral triples. In the process, we show
that for certain graphs and $k$-graphs, the spectral triples
constructed in \cite{PRen,PRS} satisfy these conditions, making them
reasonable candidates for the title of noncommutative manifolds.
Conditions for noncompact noncommutative manifolds have previously
been considered in \cite{R1,R2,GGISV}.

We have made an effort to generalise the conditions from \cite{RV}
in as minimal and stringent a way as possible. Nevertheless, our
conditions must be regarded as provisional. Additional examples and
applications are required to determine the `correct' conditions
characterising noncommutative manifolds.

The vast majority of examples of noncommutative manifolds in this
paper come from nonunital algebras (see \cite{PRen, PRS}), so our
conditions must address aspects of `noncompact noncommutative
manifolds.' Moreover, most of our examples are semifinite, in that
the trace employed is not the operator trace on Hilbert space; it is
a faithful normal semifinite trace on a different von Neumann
algebra. This is not to say that the $C^*$-algebras arising in our
examples do not admit type~I spectral triples. By considering traces
which reflect the geometry of the underlying graph (or $k$-graph),
we are naturally led to semifinite spectral triples.

For simplicity we discuss only graph algebras (i.e. algebras of
$1$-graphs) in detail; in a final section we summarise the $k$-graph
situation, since it is largely similar.

The conditions introduced in Section~\ref{sec:conditions} do not
employ Poincar\'{e} Duality in $K$-theory, but rather the Morita
equivalence condition characterising spin$^c$ structures, \cite{P}.
The equivalence of these two conditions in the compact commutative
case (in the presence of the other conditions)  was proved in
\cite{RV}. In addition, we do not consider the metric condition,
since this has recently been shown to be redundant~\cite{RV2}.

{\bf Acknowledgements} We thank J. Varilly for useful discussions.
This work was supported by the ARC and the Danish Research Council.

\section{Background on Graph $C^*$-Algebras and Spectral Triples}

\subsection{The $C^*$-algebras of Graphs}\label{graphalg}

For a more detailed introduction to graph $C^*$-algebras we refer
the reader to \cite{BPRS, kpr, CBMSbk} and the references therein. A
directed graph $E=(E^0,E^1,r,s)$ consists of countable sets $E^0$ of
vertices and $E^1$ of edges, and maps $r,s:E^1\to E^0$ identifying
the range and source of each edge. The graph is \emph{row-finite} if
each vertex emits at most finitely many edges and \emph{locally
finite} if it is row-finite and each vertex receives at most
finitely many edges. We write $E^n$ for the set of paths
$\mu=\mu_1\mu_2\cdots\mu_n$ of length $|\mu|:=n$; that is, sequences
of edges $\mu_i$ such that $r(\mu_i)=s(\mu_{i+1})$ for $1\leq i<n$.
The maps $r,s$ extend to $E^*:=\bigcup_{n\ge 0} E^n$ in an obvious
way. A \emph{loop} in $E$ is a path $L \in E^*$ with $s ( L ) = r (
L )$, we say that a loop $L$ has an exit if there is $v = s ( L_i )$
for some $i$ which emits more than one edge. If $V \subseteq E^0$
then we write $V \ge w$ if there is a path $\mu \in E^*$ with $s (
\mu ) \in V$ and $r ( \mu ) = w$. A \emph{sink} is a vertex $v \in
E^0$ with $s^{-1} (v) = \emptyset$, a \emph{source} is a vertex $w
\in E^0$ with $r^{-1} (w) = \emptyset$.

A \emph{Cuntz-Krieger $E$-family} in a $C^*$-algebra $B$ consists
of mutually orthogonal projections $\{p_v:v\in E^0\}$ and partial
isometries $\{S_e:e\in E^1\}$ satisfying the \emph{Cuntz-Krieger
relations}
\begin{equation*}
S_e^* S_e=p_{r(e)} \mbox{ for $e\in E^1$} \ \mbox{ and }\
p_v=\sum_{\{ e : s(e)=v\}} S_e S_e^* .  \mbox{ whenever $v$ is not
a sink.}
\end{equation*}

Theorem~1.2 of \cite{kpr} shows that there is a universal
$C^*$-algebra $C^*(E)$ generated by a universal Cuntz-Krieger
$E$-family $\{S_e,p_v\}$.  A product $S_\mu:=S_{\mu_1}S_{\mu_2}\dots
S_{\mu_n}$ is non-zero precisely when $\mu=\mu_1\mu_2\cdots\mu_n$ is
a path in $E^n$. Since the Cuntz-Krieger relations imply that the
projections $S_eS_e^*$ are also mutually orthogonal, we have
$S_e^*S_f=0$ unless $e=f$, and words in $\{S_e,S_f^*\}$ collapse to
products of the form $S_\mu S_\nu^*$ for $\mu,\nu\in E^*$ satisfying
$r(\mu)=r(\nu)$ (cf.\ \cite[Lemma
  1.1]{kpr}).
Indeed, because the family $\{S_\mu S_\nu^*\}$ is closed under
multiplication and involution, we have
\begin{equation}
C^*(E)=\clsp\{S_\mu S_\nu^*:\mu,\nu\in E^*\mbox{ and
}r(\mu)=r(\nu)\}.\label{spanningset}
\end{equation}
The algebraic relations and the density of $\mbox{span}\{S_\mu
S_\nu^*\}$ in $C^*(E)$ play a critical role throughout the paper.
We adopt the conventions that vertices are paths of length 0, that
$S_v:=p_v$ for $v\in E^0$, and that all paths $\mu,\nu$ appearing
in (\ref{spanningset}) are non-empty; we recover $S_\mu$, for
example, by taking $\nu=r(\mu)$, so that $S_\mu S_\nu^*=S_\mu
p_{r(\mu)}=S_\mu$.

If $z\in S^1$, then the family $\{zS_e,p_v\}$ is another
Cuntz-Krieger $E$-family which generates $C^*(E)$, and the
universal property gives a homomorphism $\gamma_z:C^*(E)\to
C^*(E)$ such that $\gamma_z(S_e)=zS_e$ and $\gamma_z(p_v)=p_v$.
The homomorphism $\gamma_{\overline z}$ is an inverse for
$\gamma_z$, so $\gamma_z\in\Aut C^*(E)$, and a routine
$\epsilon/3$ argument using (\ref{spanningset}) shows that
$\gamma$ is a strongly continuous action of $S^1$ on $C^*(E)$. It
is called the \emph{gauge action}. Because $S^1$ is compact,
averaging over $\gamma$ with respect to normalised Haar measure
gives an expectation $\Phi$ of $C^*(E)$ onto the fixed-point
algebra $C^*(E)^\gamma$:
\[
\Phi(a):=\frac{1}{2\pi}\int_{S^1} \gamma_z(a)\,d\theta\ \mbox{ for
}\ a\in C^*(E),\ \ z=e^{i\theta}.
\]
The map $\Phi$ is positive, has norm $1$, and is faithful in the
sense that $\Phi(a^*a)=0$ implies $a=0$.

From Equation (\ref{spanningset}), it is easy to see that a graph
$C^*$-algebra is unital if and only if the underlying graph is
finite. When we consider infinite graphs, formulas which involve
sums of projections may contain infinite sums. To interpret these,
we use strict convergence in the multiplier algebra of $C^*(E)$. The
following is proved in \cite{PR}.

\begin{lemma}\label{strict}
Let $E$ be a row-finite graph, let $A$ be a $C^*$-algebra
generated by a Cuntz-Krieger $E$-family $\{T_e,q_v\}$, and let
$\{p_n\}$ be a sequence of projections  in $A$. If $p_nT_\mu
T_\nu^*$ converges for every $\mu,\nu\in E^*$, then $\{p_n\}$
converges strictly to a projection $p\in M(A)$.
\end{lemma}

Another graph theoretic concept useful for the graphs we will be
dealing with is the following.
\begin{defn}\label{ends} Let $E$ be a row-finite directed graph.
 An {\em end} will mean a sink, a loop without exit or an infinite path
with no exits.
\end{defn}

Ends play two  roles in this paper. One is connected to the
existence of faithful gauge invariant traces. The most basic
statement about such traces is the following, \cite{PRen}.

\begin{lemma}\label{necessary} Let $E$ be a row-finite directed
  graph.
\par\noindent {\bf (i)} If $C^*(E)$ has a faithful, semifinite
trace then no loop
can have an exit.
\par\noindent
{\bf (ii)} If $C^* (E)$ has a faithful, semifinite, lower
semicontinuous, gauge-invariant trace $\tau$ then
$\tau \circ \Phi = \tau$ and
$$
\tau(S_\mu S_\nu^*)=\delta_{\mu,\nu}\tau(p_{r(\mu)}).
$$
\end{lemma}
This result was sharpened using the notion of a graph trace, which
we also require in this paper.

\begin{defn}[\cite{T}] If $E$ is a row-finite directed graph, then a
graph trace on $E$ is a function $g:E^0\to{\R}^+$ such that for
any $v\in E^0$ we have
\begin{equation} \label{tracecond}
g(v)=\sum_{s(e)=v}g(r(e)).
\end{equation}
\end{defn}
We showed that if $E$ admits a faithful graph trace then no vertex
$v\in E^0$ connects to any other vertex or any end via infinitely
many distinct paths. In particular, $E$ can contain no loops with
exits. We also shows that if $E$ contains finitely many ends and
every infinite path $x\in E^\infty$ is eventually contained in an
end, then $E$ admits a faithful graph trace.

The following proposition from \cite{PRen} shows how to use the
existence of graph traces on $E$ to deduce the existence of
faithful, semifinite, lower semicontinuous, gauge invariant traces
on $C^*(E)$.

\begin{prop}\label{trace=graphtrace} Let $E$ be a row-finite directed graph.
Then there is a one-to-one correspondence between faithful graph
traces on $E$ and faithful, semifinite, lower semicontinuous,
gauge invariant traces on $C^*(E)$.
\end{prop}

The other role of ends for graphs such that no loop has an exit is
topological, and will be used when we discuss the conditions. The
following is proved in \cite[Corollary 5.3]{RSz}.

\begin{lemma}\label{Kofgraph} Let $A=C^*(E)$ be a graph $C^*$-algebra
such that no
loop in the locally finite graph $E$ has an exit.  Then,
$$K_0(C^*(E))=\Z^{\#ends},\ \ \ \ K_1(C^*(E))=
\Z^{\#loops}.$$ In particular, $K_*(C^*(E))$ is finitely generated
if there are finitely many ends in $E$.
\end{lemma}

\subsection{Semifinite Spectral Triples}

We begin with some semifinite versions of standard definitions and
results. Let $\tau$ be a fixed faithful, normal, semifinite trace
on the von Neumann algebra ${\mathcal N}$. Let ${\mathcal
K}_{\mathcal N }$ be the $\tau$-compact operators in ${\mathcal
N}$ (that is the norm closed ideal generated by the projections
$E\in\mathcal N$ with $\tau(E)<\infty$).

\begin{defn} A semifinite
spectral triple $(\A,\HH,\D)$ consists of a Hilbert space $\HH$, a
$*$-algebra $\A\subset \cn$ where $\cn$ is a semifinite von Neumann
algebra acting on $\HH$, and a densely defined unbounded
self-adjoint operator $\D$ affiliated to $\cn$ such that

1) $[\D,a]$ is densely defined and extends to a bounded operator in
$\cn$ for all $a\in\A$, and

2) $a(\lambda-\D)^{-1}\in\K_\cn$ for all $\lambda\not\in{\R}$ and
all $a\in\A$.

The triple is said to be even if there is $\Gamma\in\cn$ such that
$\Gamma^*=\Gamma$, $\Gamma^2=1$,  $a\Gamma=\Gamma a$ for all
$a\in\A$ and $\D\Gamma+\Gamma\D=0$. Otherwise it is odd.
\end{defn}

\begin{defn}\label{qck} A semifinite spectral triple $(\A,\HH,\D)$ is $QC^k$ for $k\geq 1$
($Q$ for quantum) if for all $a\in\A$ the operators $a$ and
$[\D,a]$ are in the domain of $\delta^k$, where
$\delta(T)=[\dd,T]$ is the partial derivation on $\cn$ defined by
$\dd$. We say that $(\A,\HH,\D)$ is $QC^\infty$ if it is $QC^k$
for all $k\geq 1$.
\end{defn}

{\bf Note}. The notation is meant to be analogous to the classical
case, but we introduce the $Q$ so that there is no confusion
between quantum differentiability of $a\in\A$ and classical
differentiability of functions.

\noindent{\bf Remarks concerning derivations and commutators}.  By
partial derivation we mean that $\delta$ is defined on some
subalgebra of $\cn$ which need not be (weakly) dense in $\cn$.
More precisely, $\mbox{dom}\ \delta=\{T\in\cn:\delta(T)\mbox{ is
bounded}\}$. We also note that if $T\in{\mathcal N}$, one can show
that $[\dd,T]$ is bounded if and only if $[(1+\D^2)^{1/2},T]$ is
bounded, by using the functional calculus to show that
$\dd-(1+\D^2)^{1/2}$ extends to a bounded operator in $\cn$. In
fact, writing $\dd_1=(1+\D^2)^{1/2}$ and $\delta_1(T)=[\dd_1,T]$
we have \ben \mbox{dom}\ \delta^n=\mbox{dom}\ \delta_1^n\ \ \ \
\mbox{for all}\ n.\een We also observe that if $T\in\cn$ and
$[\D,T]$ is bounded, then $[\D,T]\in\cn$. Similar comments apply
to $[\dd,T]$, $[(1+\D^2)^{1/2},T]$. The proofs of these statements
can be found in \cite{CPRS2}.
\begin{defn}A $*$-algebra $\A$ is smooth if it is Fr\'{e}chet
and $*$-isomorphic to a proper dense subalgebra $i(\A)$ of a
$C^*$-algebra $A$ which is stable under the holomorphic functional
calculus.\end{defn} Thus saying that $\A$ is \emph{smooth} means
that $\A$ is Fr\'{e}chet and a pre-$C^*$-algebra. Asking for
$i(\A)$ to be a {\it proper} dense subalgebra of $A$ immediately
implies that the Fr\'{e}chet topology of $\A$ is finer than the
$C^*$-topology of $A$ (since Fr\'{e}chet means locally convex,
metrizable and complete.) We will sometimes speak of
$\overline{\A}=A$, particularly when $\A$ is represented on
Hilbert space and the norm closure $\overline{\A}$ is unambiguous.
At other times we regard $i:\A\hookrightarrow A$ as an embedding
of $\A$ in a $C^*$-algebra. We will use both points of view.

It has been shown that if $\A$ is smooth in $A$ then $M_n(\A)$ is
smooth in $M_n(A)$, \cite{GVF,LBS}. This ensures that the
$K$-theories of the two algebras are isomorphic, the isomorphism
being induced by the inclusion map $i$. This definition ensures
that a smooth algebra is a `good' algebra, \cite{GVF}, so these
algebras have a sensible spectral theory which agrees with that
defined using the $C^*$-closure, and the group of invertibles is
open.


Stability under the holomorphic functional calculus extends to
nonunital algebras, since the spectrum of an element in a
nonunital algebra is defined to be the spectrum of this element in
the  `one-point' unitization, though we must of course restrict to
functions satisfying $f(0)=0$. Likewise, the definition of a
Fr\'{e}chet algebra does not require a unit. The point of contact
between smooth algebras and $QC^\infty$ spectral triples is the
following Lemma, proved in \cite{R1}.

\begin{lemma}\label{smo} If $(\A,\HH,\D)$ is a $QC^\infty$ spectral triple, then
$(\A_\delta,\HH,\D)$ is also a $QC^\infty$ spectral triple, where
$\A_\delta$ is the completion of $\A$ in the locally convex
topology determined by the seminorms \ben
q_{n,i}(a)=\n\delta^nd^i(a)\n,\ \ n\geq 0,\ i=0,1,\een where
$d(a)=[\D,a]$. Moreover, $\A_\delta$ is a smooth algebra.
\end{lemma}

We call the topology on $\A$ determined by the seminorms $q_{ni}$
of Lemma \ref{smo} the $\delta$-topology.

Whilst smoothness does not depend on whether $\A$ is unital or
not, many analytical problems arise because of the lack of a unit.
As in \cite{R1,R2,GGISV}, we make two definitions to address these
issues.

\begin{defn} An algebra $\A$ has local units if for every finite subset of
elements $\{a_i\}_{i=1}^n\subset\A$, there exists $\phi\in\A$ such
that for each $i$ \ben \phi a_i= a_i\phi=a_i.\een
\end{defn}

\begin{defn}
Let $\A$ be a Fr\'{e}chet algebra and $\A_c\subseteq\A$ be a dense
subalgebra with local units. Then we call  $\A$ a quasi-local
algebra (when $\A_c$ is understood.) If $\A_c$ is a dense ideal
with local units, we call $\A_c\subset\A$ local.
\end{defn}

Quasi-local algebras have an approximate unit $\{\phi_n\}_{n\geq
1}\subset\A_c$ such that $\phi_{n+1}\phi_n=\phi_n$, \cite{R1}.

{\bf Example} For a graph $C^*$-algebra $A=C^*(E)$, Equation
(\ref{spanningset}) shows that
$$ A_c=\mbox{span}\{S_\mu S_\nu^*:\mu,\nu\in E^*\ \mbox{and}\
r(\mu)=r(\nu)\}$$ is a dense subalgebra. It has local units
because
$$ p_{v}S_\mu S_\nu^*=\left\{\begin{array}{ll} S_\mu S_\nu^* &
v=s(\mu)\\ 0 & \mbox{otherwise}\end{array}\right..$$ Similar
comments apply to right multiplication by $p_{s(\nu)}$. By summing
the source and range projections (without repetitions) of all
$S_{\mu_i}S_{\nu_i}^*$ appearing in a finite sum
$$ a=\sum_ic_{\mu_i,\nu_i}S_{\mu_i}S_{\nu_i}^*$$
we obtain a local unit for $a\in A_c$. By repeating this process
for any finite collection of such $a\in A_c$ we see that $A_c$ has
local units.

We also require that when we have a spectral triple the operator
$\D$ is compatible with the quasi-local structure of the algebra,
in the following sense.

\begin{defn} If $(\A,\HH,\D)$ is a spectral triple, then we define $\cda$
to be the algebra generated by $\A$ and $[\D,\A]$.
\end{defn}

\begin{defn}\label{lst} A local spectral triple $(\A,\HH,\D)$ is a
spectral triple with $\A$ quasi-local such that there exists an
approximate unit $\{\phi_n\}\subset\A_c$ for $\A$ satisfying \ben
{\mathcal C}_\D(\A_c)=\bigcup_n\cda_n,\een where \ben
\cda_n=\{\omega\in\cda:\phi_n\omega=\omega\phi_n=\omega\}.\een
\end{defn}

{\bf Remark} A local spectral triple has a local approximate unit
$\{\phi_n\}_{n\geq 1}\subset\A_c$  such that
$\phi_{n+1}\phi_n=\phi_n\phi_{n+1}=\phi_n$ and
$\phi_{n+1}[\D,\phi_n]=[\D,\phi_n]\phi_{n+1}=[\D,\phi_n]$,
\cite{R2}. This is the crucial property we require to prove most of
our summability results, to which we now turn.

\subsection{Summability}
In the following, let $\mathcal N$ be a semifinite von Neumann
algebra with faithful normal trace $\tau$. Recall from \cite{FK}
that if $S\in\mathcal N$, the \emph{t-th generalized singular value}
of S for each real $t>0$ is given by
$$\mu_t(S)=\inf\{||SE||\ : \ E \mbox{ is a projection in }
{\mathcal N} \mbox { with } \tau(1-E)\leq t\}.$$

The ideal $\LL^1({\mathcal N})$ consists of those operators $T\in
{\mathcal N}$ such that $\|T\|_1:=\tau( |T|)<\infty$ where
$|T|=\sqrt{T^*T}$. In the Type I setting this is the usual trace
class ideal. We will simply write $\LL^1$ for this ideal in order to
simplify the notation, and denote the norm on $\LL^1$ by
$\n\cdot\n_1$. An alternative definition in terms of singular values
is that $T\in\LL^1$ if $\|T\|_1:=\int_0^\infty \mu_t(T) dt <\infty.$

Note that in the case where ${\mathcal N}\neq{\mathcal
B}({\mathcal H})$, $\LL^1$ is not complete in this norm but it is
complete in the norm $||.||_1 + ||.||_\infty$. (where
$||.||_\infty$ is the uniform norm). Another important ideal for
us is the domain of the Dixmier trace:
$${\mathcal L}^{(1,\infty)}({\mathcal N})=
\left\{T\in{\mathcal N}\ : \Vert T\Vert_{_{{\mathcal
L}^{(1,\infty)}}} :=   \sup_{t> 0}
\frac{1}{\log(1+t)}\int_0^t\mu_s(T)ds<\infty\right\}.$$


We will suppress the $({\mathcal N})$ in our notation for these
ideals, as $\cn$ will always be clear from context. The reader
should note that ${\mathcal L}^{(1,\infty)}$ is often taken to
mean an ideal in the algebra $\widetilde{\mathcal N}$ of
$\tau$-measurable operators affiliated to ${\mathcal N}$. Our
notation is however consistent with that of \cite{C} in the
special case ${\mathcal N}={\mathcal B}({\mathcal H})$. With this
convention the ideal of $\tau$-compact operators, ${\mathcal
  K}({\mathcal N})$,
consists of those $T\in{\mathcal N}$ (as opposed to
$\widetilde{\mathcal N}$) such that \ben \mu_\infty(T):=\lim
_{t\to \infty}\mu_t(T)  = 0.\een

\begin{defn}\label{summable} A semifinite local spectral triple is
$(1,\infty)$-summable if \ben
a(\D-\lambda)^{-1}\in\LL^{(1,\infty)}\ \ \ \mbox{for all}\
a\in\A_c,\ \ \lambda\in\C\setminus\R.\een
\end{defn}

{\bf Remark} If $\A$ is unital, $\ker\D$ is $\tau$-finite
dimensional. Note that the summability requirements are only for
$a\in\A_c$. We do not assume that elements of the algebra $\A$ are
all integrable in the nonunital case. Strictly speaking, this
definition describes {\em local} $(1,\infty)$-summability, however
we use the terminology $(1,\infty)$-summable to be consistent with
the unital case.

We need to briefly discuss the Dixmier trace, but fortunately we
will usually be applying it in reasonably simple situations. For
more information on semifinite Dixmier traces, see \cite{CPS2}.
For $T\in\LL^{(1,\infty)}$, $T\geq 0$, the function \ben
F_T:t\to\frac{1}{\log(1+t)}\int_0^t\mu_s(T)ds \een is bounded. For
certain generalised limits $\omega\in L^\infty(\R_*^+)^*$, we
obtain a positive functional on $\LL^{(1,\infty)}$ by setting
$$ \tau_\omega(T)=\omega(F_T).$$
 This is the
Dixmier trace associated to the semifinite normal trace $\tau$,
denoted $\tau_\omega$, and we extend it to all of
$\LL^{(1,\infty)}$ by linearity, where of course it is a trace.
The Dixmier trace $\tau_\omega$ is defined on the ideal
$\LL^{(1,\infty)}$, and vanishes on the ideal of trace class
operators. Whenever the function $F_T$ has a limit at infinity,
all Dixmier traces return the value of the limit. We denote the
common value of all Dixmier traces on measurable operators by
$\bigintcross$. So if $T\in\LL^{(1,\infty)}$ is measurable, for
any allowed functional $\omega\in L^\infty(\R_*^+)^*$ we have
$$\tau_\omega(T)=\omega(F_T)=\bigintcross T.$$

{\bf Example} Let $\D=\frac{1}{i}\frac{d}{d\theta}$ act on
$L^2(S^1)$. Then it is well known that the spectrum of $\D$ consists
of eigenvalues $\{n\in\Z\}$, each with multiplicity one. So, using
the standard operator trace, $\mbox{Trace}$, the function
$F_{(1+\D^2)^{-1/2}}$ is
$$ \frac{1}{\log 2N}\sum_{n=-N}^N(1+n^2)^{-1/2}$$
and this is bounded. Hence $(1+\D^2)^{-1/2}\in\LL^{(1,\infty)}$ and
\begin{equation}\mbox{Trace}_\omega((1+\D^2)^{-1/2})=\bigintcross(1+\D^2)^{-1/2}=2.\label{eq:2}\end{equation}

In \cite{R1,R2} we proved numerous properties of local algebras.
The introduction of quasi-local algebras in \cite{GGISV} led us to
review the validity of many of these results for quasi-local
algebras. Most of the summability results of \cite{R2} are valid
in the quasi-local setting.  In addition, the summability results
of \cite{R2} are also valid for general semifinite spectral
triples since they rely only on properties of the ideals
$\LL^{(p,\infty)}$, $p\geq 1$, \cite{C,CPS2}, and the trace
property. We quote the version of the summability results from
\cite{R2} that we require below.

\begin{prop}[\cite{R2}]\label{wellbehaved} Let $(\A,\HH,\D)$ be a $QC^\infty$, local
$(1,\infty)$-summable semifinite spectral triple relative to
$(\cn,\tau)$. Let $T\in\cn$ satisfy $T\phi=\phi T=T$ for some
$\phi\in\A_c$. Then \ben T(1+\D^2)^{-1/2}\in\LL^{(1,\infty)}.\een
For $Re(s)>1$, $T(1+\D^2)^{-s/2}$ is trace class. If the limit \be
\lim_{s\to 1/2^+}(s-1/2)\tau(T(1+\D^2)^{-s})\label{mumbo}\ee
exists, then it is equal to \ben \frac{1}{2}\bigintcross
T(1+\D^2)^{-1/2}.\een In addition, for any Dixmier trace
$\tau_\omega$, the function \ben a\mapsto
\tau_\omega(a(1+\D^2)^{-1/2})\een defines a trace on
$\A_c\subset\A$.
\end{prop}

\subsection{The Gauge Spectral Triple for a Graph $C^*$-Algebra}
In this section we summarise the construction of a Kasparov module
and a semifinite spectral triple for locally finite directed graphs
with no sources. This material is based on \cite{PRen}. We begin by
constructing a Kasparov module.

For $E$ a row finite directed graph, we set $A=C^*(E)$,
$F=C^*(E)^\gamma$, the fixed point algebra for the $S^1$ gauge
action. The algebras $A_c,F_c$ are defined as the finite linear
span of the generators.

Right multiplication makes $A$ into a right $F$-module, and
similarly $A_c$ is a right module over $F_c$. We define an
$F$-valued inner product $(\cdot|\cdot)_R$ on both these modules
by
$$ (a|b)_R:=\Phi(a^*b),$$
where $\Phi$ is the canonical expectation $A\to F$.
\begin{defn}\label{mod} Let $X$ be the right $C^*$-$F$-module obtained by
completing $A$ (or $A_c$) in the norm
$$ \Vert x\Vert^2_X:=\Vert (x|x)_R\Vert_F=\Vert
\Phi(x^*x)\Vert_F.$$
\end{defn}

The algebra $A$ acting by multiplication on the left of $X$
provides a representation of $A$ as adjointable operators on $X$.
We let $X_c$ be the copy of $A_c\subset X$.

For each $k\in\Z$, define an operator $\Phi_k$ on $X$  by
$$\Phi_k(x)=\frac{1}{2\pi}\int_{S^1}z^{-k}\gamma_z(x)d\theta,\ \ z=e^{i\theta},\ \ x\in X.$$
Observe that on generators we have \begin{equation}\Phi_k(S_\al
S_\beta^*)=\left\{\begin{array}{lr}S_\al S_\beta^* & \ \
|\al|-|\beta|=k\\0 & \ \ |\al|-|\beta|\neq
k\end{array}\right..\label{kthproj}\end{equation}

\begin{lemma}\label{phiendo} The operators $\Phi_k$ are adjointable endomorphisms
of the $F$-module $X$ such that $\Phi_k^*=\Phi_k=\Phi_k^2$ and
$\Phi_k\Phi_l=\delta_{k,l}\Phi_k$. If $K\subset\Z$ then the sum
$\sum_{k\in K}\Phi_k$ converges strictly to a projection in the
endomorphism algebra. The sum $\sum_{k\in\Z}\Phi_k$ converges to
the identity operator on $X$.
\end{lemma}

\begin{cor}\label{gradedsum} Let $x\in X$. Then with $x_k=\Phi_kx$ the sum
$\sum_{k\in \Z}x_k$ converges in $X$ to $x$.
\end{cor}

\begin{prop}[\cite{PRen}] Let $X$ be the right  $C^*$-$F$-module of
Definition \ref{mod}. Let
$$ X_\D=\{x\in X:\sum_{k\in\Z}k^2(x_k|x_k)<\infty\},$$
and define $\D:X_\D\to X$ by
$$ \D\sum_{k\in\Z}x_k=\sum_{k\in\Z}kx_k.$$
Then $\D$ is closed,  self-adjoint and regular.
\end{prop}

We refer to Lance's book, \cite[Chapters 9,10]{L}, for information
on unbounded operators on $C^*$-modules.
\begin{lemma}\label{finrank} Assume that the directed graph $E$ is
locally finite and has no sources. For all $a\in A$ and $k\in\Z$,
$a\Phi_k\in End^0_F(X)$, the
 compact right endomorphisms of $X$. If $a\in A_c$ then
$a\Phi_k$ is finite rank.
\end{lemma}

\begin{thm} Suppose that the graph $E$ is locally finite and has no sources,
and let $X$ be the right $F$ module of Definition \ref{mod}. Let
$V=\D(1+\D^2)^{-1/2}$. Then $(X,V)$ is an odd Kasparov module for
$A$-$F$ and so defines an element of $KK^1(A,F)$.
\end{thm}

Given the hypotheses of the Theorem, we may describe $\D$ as
$$\D=\sum_{k\in\Z}k\Phi_k.$$

To construct a semifinite spectral triple, we suppose that our graph
$C^*$-algebra also has a faithful gauge invariant trace
$\tau:A\to\C$. Using the trace $\tau$, we define a $\C$-valued inner
product $\la\cdot,\cdot\ra$ on $X_c$ by
$$\la x,y\ra:=\tau((x|y)_R)=\tau(\Phi(x^*y))=\tau(x^*y),$$
the last equality following from the gauge invariance of $\tau$.
Denote the Hilbert space completion of $X_c$ by $\HH=L^2(X,\tau)$.

\begin{prop} The actions of $A$ by left and right multiplication
on $X$ extend to give commuting bounded faithful nondegenerate
representations of $A$ and $A^{op}$, where $A^{op}$ is the
opposite algebra of $A$. Furthermore, any endomorphism of the
right $F$ module $X$ which preserves $X_c$ extends uniquely
to a bounded operator on $\HH$.
\end{prop}

The operator $\D:\mbox{dom}\D\subset X\to X$ also extends to a
self-adjoint operator on $\HH$, \cite[Lemma 5.5]{PRen}. It is
shown in \cite{PRen} that for all $a\in A_c$ the commutator
$[\D,a]$ extends to a bounded operator on $\HH$.

\begin{lemma}\label{smoothquasi} The algebra $A_c$ and the
linear space $[\D,A_c]$ are
 contained in the smooth domain
of the derivation $\delta$ where for $T\in\B(\HH)$,
$\delta(T)=[\dd,T]$. So the completion of $A_c$ in the
$\delta$-topology, which we denote by $\A$, is a Fr\'{e}chet
pre-$C^*$-algebra. Moreover $\A$ is a quasi-local algebra with
dense subalgebra $A_c$.
\end{lemma}

Since $[\D,A_c]\subset A_c$, it is not hard to see that if
$(\A,\HH,\D)$ is a spectral triple, then it is local.

The only remaining piece of information we require to obtain a
spectral triple is the von Neumann algebra and trace which give us
a semifinite spectral triple. Let $End^{00}_F(X_c)$ be the finite rank
endomorphisms of the pre-$C^*$-module $X_c$.
\begin{prop} Let $\cn=(End^{00}_F(X_c))''$. Then there exists a
  faithful, normal, semifinite trace $\tilde\tau:\cn\to\C$ such that
  for all rank one endomorphisms $\Theta_{x,y}$ of $X_c$ we have
$$ \tilde\tau(\Theta_{x,y})=\tau((y|x)_R),\ \ x,y\in X_c.$$
Moreover, for all $a\in\A$ and $\lambda\in\C\setminus\R$ the operator
$a(\lambda-\D)^{-1}$  lies in $\K_\cn$.
\end{prop}
Hence we obtain a
semifinite spectral triple. However, more is true.

\begin{thm}\label{mainthm} Let $E$ be a locally finite graph with no sources,
and let $\tau$ be a faithful, semifinite, norm lower-semicontinuous,
gauge invariant trace on $C^*(E)$. Then $(\A,\HH,\D)$ is a
$QC^\infty$ $(1,\infty)$-summable odd local semifinite spectral
triple (relative to $(\cn,\tilde\tau)$). For all $a\in \A$ the
operator $a(1+\D^2)^{-1/2}$ is not trace class. For any $v\in E^0$
which does not connect to a sink we have
$$\tilde\tau_\omega(p_v(1+\D^2)^{-1/2})=2\tau(p_v),$$
where $\tilde\tau_\omega$ is any Dixmier trace associated to
$\tilde\tau$.
\end{thm}

The main observation of the proof is that for $v\in E^0$ such that
$v$ does not connect to a sink, and for $k\in\Z$ we have
$$\tilde\tau(p_v\Phi_k)=\tau(p_v).$$

This is the spectral triple we will be working with for the rest of
the paper, and we refer to it as the gauge spectral triple of the
directed graph $E$ (or algebra $C^*(E)$). We remind the reader that
the existence of this spectral triple depends only on the graph $E$
being locally finite with no sources, and the existence of a
faithful, semifinite, gauge invariant, norm lower-semicontinuous
trace $\tau:A\to\C$. The latter is a nontrivial condition.
\section{Conditions for Locally Compact Semifinite
Manifolds}\label{sec:conditions}

We now review in turn the conditions for noncommutative manifolds as
presented in \cite{RV}. We will consider natural generalisations of
these conditions to the semifinite and nonunital setting and
consider when the gauge spectral triple satisfies these
generalisations.

We will present each condition as stated for the type I and unital
case, where $(\cn,\tau)=(\B(\HH),\mbox{Trace})$ and
$(1+\D^2)^{-1/2}\in\K(\HH)$, and then present the necessary
modification of the condition, if it requires modification.

When dealing with these generalisations we will suppose that
$(\A,\HH,\D)$ is a local semifinite spectral triple relative to
$\cn,\tilde\tau$. We will not suppose that $\A$ is unital, but will
suppose that $\A_c\subset\A$ gives us a quasi-local algebra.

When considering the conditions as applied to graph algebras, we
will suppose that  $E$ is a locally finite directed graph with no
sources and possessing a faithful graph trace $g$. We will let
$(\A,\HH,\D)$ be the gauge spectral triple of $E$ described in the
previous section.

The conditions are somewhat interdependent, and we have found it is
difficult to present them in a logical fashion. It seems that this
difficulty is greatly eased if we assume at the outset that the
Hilbert space $\HH$ carries commuting representations
$\pi:\A\to\B(\HH)$ and $\pi^{op}:\A^{op}\to\B(\HH)$. The former
representation actually has $\pi(\A)\subset\cn\subset\B(\HH)$, but
we do not assume this for the latter representation.

We will explicitly state this bimodule requirement again when we
look at the first order condition, but it will be apparent that
several of our conditions require a bimodule structure for their
statement. In all the following, we identify $a \in \A$ with $\pi(a)
\in \cn$ unless stated otherwise.

\subsection{The Analytic Conditions}

\begin{ocond}[Dimension]
\label{cn:metr-dim} The type I unital spectral triple $(\A,\HH,\D)$
is $(p,\infty)$-summable for a fixed positive \textit{integer} $p$,
for which ${\rm Trace}_\omega((1 + D^2)^{-p/2}) > 0$ for all Dixmier
limits $\omega$.
\end{ocond}

To generalise  this condition we evidently need to replace the
operator trace, $\mbox{Trace}$, by the trace $\tilde\tau:\cn\to\C$
which determines the compactness and summability requirements of our
spectral triple. We also need to restate the requirement, since in
general for a nonunital spectral triple, even type I, we will not
have $(1+\D^2)^{-p/2}\in\LL^{(1,\infty)}$, \cite{R2}. So we have a
simultaneous generalisation to the nonunital and semifinite case.

\begin{ncond}[Semifinite Nonunital Dimension]
\label{cn:metr-dim-sf} The local semifinite spectral triple
$(\A,\HH,\D)$ is $(p,\infty)$-summable for a fixed positive
\textit{integer} $p$, for which ${\rm \tilde\tau}_\omega(a(1 +
D^2)^{-p/2}) > 0$ for all $\omega$ and all $a\in\A_c$ with $a\geq
0$.
\end{ncond}

This generalisation of the dimension condition is satisfied by the
gauge spectral triple of a directed graph
with $p=1$. Provided the graph $E$ has no
sinks this follows from Theorem \ref{mainthm} since the Dixmier
trace of $a^*a$, $0\neq a\in\A_c$, is given by
$\tilde\tau_\omega(a^*a(1+\D^2)^{-1/2})=2\tau(a^*a)>0$. Even if
the graph has sinks, the proof of Theorem \ref{mainthm} in
\cite{PRen} shows that we still have positivity.

\begin{ocond}[Regularity]
\label{cn:qc-infty} The spectral triple $(\A,\HH,\D)$ is
$QC^\infty$. Without loss of generality, we assume that $\A$ is
complete in the $\delta$-topology and so is a Fr\'echet
pre-$C^*$-algebra.
\end{ocond}

It follows from Lemma \ref{smoothquasi} that this condition is
satisfied with no need to modify it at all.

\begin{ncond}[Regularity]
The spectral triple $(\A,\HH,\D)$ is $QC^\infty$. Without loss of
generality, we assume that $\A$ is complete in the $\delta$-topology
and so is a Fr\'echet pre-$C^*$-algebra.
\end{ncond}

\textbf{Remark.} In the type I setting, we also have the condition
of \emph{Absolute Continuity} which states: \emph{For all nonzero $a
\in \A$ with $a \geq 0$, and for any $\omega$-limit, the following
Dixmier trace is positive:
$$
{\rm Trace}_\omega(a(1 + \D^2)^{-p/2}) > 0.
$$}

This is half of Connes' finiteness and absolute continuity
condition, \cite{C1,GVF}, the other half being finiteness discussed
in Section~\ref{sec:bimodule} below; see also \cite{RV}. Here we
have demanded positivity only for positive elements of $\A_c$, but
this extends to positive elements of $\A$, provided we allow the
value $+\infty$. Of course our reformulation of the dimension
condition already subsumes a semifinite version of absolute
continuity, so the natural generalisation of the absolute continuity
condition is already satisfied by our gauge spectral triples. This
shows that even in the unital case it makes sense to combine the
dimension and absolute continuity conditions, as mentioned in
\cite{RV}.

Thus our formulation of the conditions has rendered the absolute
continuity condition redundant.

\subsection{The Orientation and Closedness Conditions}

This section examines the orientation and finiteness conditions. The
orientability condition for spectral triples with unital algebra
$\A$ is

\begin{ocond}[Orientability]
\label{cn:orient} Let $p$ be the metric dimension of $(\A,\HH,\D)$.
We require that the spectral triple be \textit{even}, with
$\Z_2$-grading $\Gamma$, if and only if $p$ is even. For
convenience, we take $\Gamma = Id_\HH$ when $p$ is odd. We say the
spectral triple $(\A,\HH,\D)$ is \textit{orientable} if there exists
a Hochschild $p$-cycle
\begin{subequations}
\label{eq:Hoch-cycle}
\begin{equation}
\label{eq:Hoch-cycle-defn} c = \sum_{\al=1}^n a_\al^0\otimes
b^{op}_\al \otimes a_\al^1 \otimes\cdots\otimes a_\al^p \in
Z_p(\A,\A\otimes\A^{op})
\end{equation}
whose Hochschild class $[c] \in HH_p(\A)\otimes\A^{op}$ may be
called the ``orientation'' of $(\A,\HH,\D)$, such that
\begin{equation}
\label{eq:Hoch-cycle-repn} \pi_D(c) := \sum_\al a_\al^0b^{op}_\al
\,[\D,a_\al^1] \dots [\D,a_\al^p] = \Gamma.
\end{equation}
\end{subequations}
\end{ocond}
Here  $\A\otimes\A^{op}$ is a bimodule for $\A$ via
$$ a\cdot(x\otimes y^{op})=ax\otimes y^{op},\ \ \ (x\otimes
y^{op})\cdot a=xa\otimes y^{op},\ \ \ a,x,y\in\A.$$

Now, typically, we have a nonunital algebra, and require a different
formulation. We adopt the attitude that we should have a locally
finite but possibly infinite cycle, as would be the case for a
volume form on a noncompact manifold.

\begin{ncond}[Nonunital Orientability]
\label{cn:orient-nonunital} Let $p$ be the metric dimension of
$(\A,\HH,\D)$. We require that the spectral triple be \textit{even},
with $\Z_2$-grading $\Gamma$, if and only if $p$ is even. For
convenience, we take $\Gamma = Id_\HH$ when $p$ is odd. We say the
spectral triple $(\A,\HH,\D)$ is \textit{orientable} if there exists
a Hochschild $p$-cycle
\begin{subequations}
\label{eq:Hoch-cycle-nonunital}
\begin{equation}
\label{eq:Hoch-cycle-defn-nonunital} c = \sum_{\al=1}^\infty
a_\al^0\otimes b^{op}_\al \otimes a_\al^1 \otimes\cdots\otimes
a_\al^p
\end{equation}
whose Hochschild class $[c]$ may be called the ``orientation'' of
$(\A,\HH,\D)$, such that
\begin{equation}
\label{eq:Hoch-cycle-repn-nonunital} \pi_D(c):= \sum_\al
a_\al^0b^{op}_\al \,[\D,a_\al^1] \dots [\D,a_\al^p] = \Gamma
\end{equation}
where the sum in~\eqref{eq:Hoch-cycle-repn-nonunital} converges
strongly.
\end{subequations}
\end{ncond}
{\bf Remark} We have deliberately omitted any mention of the
homology groups that $c$ should belong to, there being many
possibilities and few examples to guide us. We offer one possible
candidate, without examining the subject in detail.

Let $C_n(A_c,A_c\otimes A_c^{op})$ be the linear space of algebraic
Hochschild $n$-chains for $A_c$. Suppose $\A$ is the completion of
$A_c$ in the topology determined by the seminorms $q_k$, let
$\{q_{{\bf k}}\}_{{\bf k}\in \N^n}$ be a corresponding family of
seminorms on $C_n(A_c,A_c\otimes A_c^{op})$ and let $\{\phi_j\}$ be
a local approximate unit for $\A$, \cite{R1}. Define
$C_n^\infty(\A,\A\otimes\A^{op})$ to be the completion of
$C_n(A_c,A_c\otimes A_c^{op})$ for the topology determined by the
family of seminorms
$$ q_{{\bf k},j}(a^0\otimes b^{op}\otimes a^1\otimes\cdots\otimes
a^n):=q_{{\bf k}}(\phi_ja^0\otimes(\phi_jb)^{op}\otimes
\phi_ja^1\otimes\cdots\otimes\phi_ja^n).$$ This should be viewed as
similar to uniform convergence of all derivatives on compacta, and
so analogous to a $C^\infty$ topology. Ultimately more nonunital
examples are required to clarify this issue; for more comments see
\cite{GGISV,R1,R2}. We leave these homological questions for future
investigation.

For the case of graph algebras, we consider the sum over all edges
in the graph
\begin{equation}\label{eq:graph Hochschild}
c=\sum_{e\in E^1} S_e^*\otimes S_e.
\end{equation}
Before worrying
about the convergence of this sum (in the multiplier algebra), we
apply the Hochschild boundary $b$ to find
$$ b(c)=\sum_e(S_e^*S_e-S_e S_e^*)=\sum_ep_{r(e)}-
\sum_{v\ not\ sink}p_v,$$ where we have used the Cuntz-Krieger
relation to obtain the second sum on the right-hand side. Thus if
there are no sinks, the second sum on the right-hand side converges
to the identity (in the multiplier algebra or the `one-point'
unitization).

The first sum on the right-hand side contains each vertex projection
$p_v$ with multiplicity equal to the number of edges entering it,
which we denote by $|v|_1$. Thus
$$b(c)=\sum_{v\in E^0,\ v\ not\ sink}(|v|_1p_v-p_v)+\sum_{v\ a\
sink}|v|_1p_v.$$ In particular, if each vertex has precisely one
edge entering it, and no vertex is a sink, $b(c)=0$. We say that
such a graph $E$ has no sinks, and satisfies the {\em single entry
condition}.

Observe that the single entry condition (together with the
requirement that no loop has an exit) rules out loops except for the
case where the (connected) graph comprises a single loop. The
$C^*$-algebra of a graph consisting of a simple loop on $n$ vertices
is isomorphic to $M_n(C(S^1))$. For a one-edge loop, the Hochschild
cycle $c$ is $z^{-1}\otimes z$, the usual volume form for the
circle. The single entry condition also rules out sources, so unless
our (connected) graph comprises a single loop, it is an infinite
directed tree with no sources or sinks, in which case the
$C^*$-algebra is AF \cite{kpr}.

If $E$ satisfies the single entry condition then we claim that $\sum
S_e^* \otimes S_e$ converges to a partial isometry in the multiplier
algebra of $C^* (E) \otimes C^*(E)$. Let $X_e = S_e^* \otimes S_e$
then it is clear that $X_e$ is a partial isometry in $C^* (E)
\otimes C^* (E)$ with \bean X_e X_e^* &=& ( S_e^* \otimes S_e ) (
S_e \otimes S_e^* ) = P_{r(e)} \otimes S_e S_e^* \nno X_e^* X_e &=&
( S_e \otimes S_e^* ) ( S_e^* \otimes S_e ) =  S_e S_e^* \otimes
P_{r(e)} . \eean By the relations in $C^* (E)$ the $S_e S_e^*$ are
mutually orthogonal, and then by the single entry hypothesis the
$P_{r(e)}$ are too. Hence the $X_e$ have mutually orthogonal ranges,
and a standard argument (see [PR2, Lemma 1.1] or [BPRS, Lemma 1.1])
finishes off the claim.

Using the single-entry condition, we see that the Hochschild cycle
defined in~\eqref{eq:graph Hochschild} is represented by
\begin{equation}\label{eq:sums}
\pi_\D(c)=\sum_eS_e^*[\D,S_e]=\sum_eS^*_eS_e=\sum_ep_{r(e)}=Id_\HH,
\end{equation}
showing that the new condition of orientation is satisfied for this
cycle. The sums in~\eqref{eq:sums} converge in the strict topology
as an operator on the $C^*$-module $X$, and also converge strongly
on $\HH$.

It may well be possible that there is a Hochschild cycle for a more
general family of graphs, and we are not claiming that the above
conditions are necessary for the orientability condition to hold,
only sufficient.

{\bf From now on we suppose that $E$ has no sinks and satisfies the
single entry condition}. As noted above, it follows that the algebra
$C^*(E)$ is then AF unless it is $M_n(C(S^1))$. In the AF case, $E$
is a directed tree. We record the following Lemma describing the
fixed-point algebra of the directed tree examples.

\begin{lemma}\label{describeF} Suppose that $E$ is a directed tree with no sinks
satisfying the single entry condition and having finitely many ends.
Then $F$ is an abelian algebra, isomorphic to the continuous
functions on the infinite path space $E^\infty$ of $E$. Letting $N$
denote the number of ends, each $f\in F_c$ can be represented as
$$ f=\sum_{v\in E^0}\sum_{n=1}^Nc_{v,n}S_{v,n}S_{v,n}^*,\ \ \ c_{v,n}\in\C$$
where $(v,n)$ denotes a path with source $v$ and range in the $n$-th
tail. The $C^*$-norm of such an $f$ is $\Vert
f\Vert_F^2=\sup|c_{v,n}|^2$.
\end{lemma}

\begin{proof} The assertion that $F\cong C_0(E^\infty)$ follows
from \cite{KPRR}. To see that it is possible to write $f\in F_c$ in
the above form, consider a path $\al$ with range $r(\al)$ a vertex
emitting one edge $e$. Then
$$ S_{\al e}S_{\al e}^*=S_\al S_eS_e^*S_\al^*=S_\al
p_{s(e)}S_\al^*=S_\al p_{r(\al)}S_\al^*=S_\al S_\al^*.$$ So any
$S_\al S_\al^*$ is equal to $S_\beta S_\beta^*$ where $\beta$ is an
extension of $\al$ not passing through a vertex emitting more than
one edge. If $\al$ is a path with range a vertex emitting, say, $k$
edges, $e_1,\dots,e_k$, then
$$S_\al S_\al^*=S_\al
p_{r(\al)}S_\al^*=\sum_{i=1}^kS_\al S_{e_i}S_{e_i}^*S_\al^*,$$ and
this can be subsequently extended until the next vertices emitting
more than one edge. This process terminates after finitely many
steps because there are finitely many ends. The $S_{v,n}S_{v,n}^*$
are mutually orthogonal, so
$$f^*f=\sum_v\sum_n|c_{v,n}|^2S_{v,n}S_{v,n}^*,$$
and $\Vert f\Vert_F^2=\sup|c_{v,n}|^2$.
\end{proof}

The next condition is closedness, which, in its original form, is
basically Stoke's theorem for the Dixmier trace applied to elements
of $\A\otimes\A^{op}$. The original formulation for
$(p,\infty)$-summable triples using the operator trace
$\mbox{Trace}$ is

\begin{ocond}[Closedness]
\label{cn:closed} The $(p,\infty)$-summable spectral triple
$(\A,\HH,\D)$ is closed if for any $a_1,\dots,a_p \in
\A\otimes\A^{op}$, the operator $\Gamma\,[\D,a_1]\cdots[\D,a_p](1 +
\D^2)^{-p/2}$ has vanishing Dixmier trace; thus, for any Dixmier
trace ${\rm Trace}_\omega$,
\begin{equation}
{\rm Trace}_\omega\bigl(\Gamma\,[\D,a_1]\cdots[\D,a_p](1 +
\D^2)^{-p/2}\bigr) = 0. \label{eq:closed-chain}
\end{equation}
\end{ocond}

{\bf Remark} By setting $\phi(a_0,\dots,a_p) := {\rm
Trace}_\omega\bigl(\Gamma\,a_0\,[\D,a_1]\,\cdots[\D,a_p](1 +
\D^2)^{-p/2}\bigr)$, the equation \eqref{eq:closed-chain} may be
rewritten \cite[VI.2]{C} as $B_0\phi = 0$, where $B_0$ is defined on
$(k+1)$-linear functionals by
$$
(B_0\phi)(a_1,\dots,a_k) := \phi(1,a_1,\dots,a_k) + (-1)^k
\phi(a_1,\dots,a_k,1).
$$

To see the utility of this condition, we introduce some notation so
that we can quote Lemma 3 of \cite[VI.4.$\gamma$]{C}. Let
$\Omega^*(\A)$ be the universal differential algebra of $\A$,
\cite[II.1.$\al$]{C}. Then  $\pi_\D : \Omega^*(\A) \to \cda$ defined
by $\pi_\D(a_0\delta a_1 \dots \delta a_n) =
a_0[\D,a_1]\cdots[\D,a_n]$ is a $*$-algebra representation. Denote
by $\Lambda^*_\D(\A)$ the graded differential algebra we obtain by
quotienting $\cda$ by the differential ideal
$\pi_\D(\delta(\ker\pi_\D))$, where $\delta$ is the universal
derivation on $\Omega^*(\A)$. We denote by $d$ the derivation on
$\Lambda^*_\D(\A)$. See \cite[Chap VI]{C} for more information.
Finally, let $Z^k(\A,\A^*)$ denote the Hochschild cocycles.

\begin{lemma}
\label{lm:PDin-cohom} Let $(\A,\HH,\D)$ be $(p,\infty)$-summable and
satisfy Old Condition~\ref{cn:first-ord} (first order). Then for
each $k = 0,1,\dots,p$ and $\eta \in \Omega^k\A$, a Hochschild
cocycle $C_\eta \in Z^{p-k}(\A,\A^*)$ is defined by
$$
C_\eta(a^0,\dots,a^{p-k}) := {\rm Trace}_\omega\bigl(
\Gamma\,\pi_\D(\eta)\, a^0\,[\D,a^1]\dots [\D,a^{p-k}]\,(1 +
\D^2)^{-p/2} \bigr).
$$
Moreover, if Old Condition~\ref{cn:closed} (closedness) also holds,
then $C_\eta$ depends only on the class of $\pi_\D(\eta)$ in
$\Lambda_\D^k \A$, and
$$
B_0 C_\eta = (-1)^k \,C_{d\eta}. \eqno\qed
$$
\end{lemma}
Thus the first order condition together with closedness give us
tools to study the Hochschild and cyclic homology of the algebra
$\A$. More information can be found in \cite[VI.4.$\gamma$]{C}.

The difficulty we face is that we have a Dixmier trace defined on
$\cn\supset\A$ which we can not apply to $\A\otimes\A^{op}$.
Nevertheless, as we discuss further in the next section, we do not
believe having a spectral triple for $\A\otimes\A^{op}$ is of
central importance. Nevertheless, we see that the utility of Lemma
\ref{lm:PDin-cohom} is greatly reduced by our new formulation.

\begin{ncond}[Semifinite Closedness]
\label{cn:sfclosed} The $(p,\infty)$-summable local semifinite
spectral triple $(\A,\HH,\D)$ is closed if for any Dixmier trace
$\tilde\tau_\omega$ we have
\begin{equation}
\tilde\tau_\omega\bigl(\Gamma\,[\D,a_1]\cdots[\D,a_p](1 +
\D^2)^{-p/2}\bigr) = 0 \label{eq:closed-chain-sf}
\end{equation}
for all $a_1,\dots,a_p \in \A$.
\end{ncond}
It would seem that this formulation does not give us tools to study
the Hochschild and cyclic cohomology of $\A$ as in the type I case
described above, \cite[VI.4.$\gamma$]{C}.

Returning to the gauge spectral triple of a graph algebra, for a
generator $S_\mu S_\nu^*\in\A$, we have
\begin{align*} \tilde\tau_\omega([\D,S_\mu S_\nu^*](1+\D^2)^{-1/2})
&=(|\mu|-|\nu|)\tilde\tau(S_\mu S_\nu^*(1+\D^2)^{-1/2})\nno
&=(|\mu|-|\nu|)\tau(S_\mu S_\nu^*).
\end{align*}
The gauge invariance of the trace says that $\tau(S_\mu S_\nu^*)$ is
non-zero only if $|\mu|=|\nu|$, whence the whole expression always
vanishes. Hence the new closedness condition holds for the gauge
spectral triple of graph algebras.

\subsection{The Bimodule Conditions}\label{sec:bimodule}
This section is concerned with the relation between the bimodule
structure of the Hilbert space and the spectral triple.

The condition of finiteness in the unital case is
\begin{ocond}[Finiteness]
\label{cn:finite} The dense subspace of $\HH$ which is the smooth
domain of~$\D$,
$$
\HH_\infty := \bigcap_{m\geq 1} {\rm dom}\ \D^m
$$
is a finitely generated projective right $\A$-module.
\end{ocond}
Thus $\HH_\infty \simeq q\,\A^m$ where $q \in M_m(\A)$ is an
idempotent. Without loss of generality, we may suppose $q = q^*$
also, so that, without further hypotheses, $\HH_\infty$ carries an
$\A$-valued Hermitian pairing, namely, that given by $(\xi,\eta)' :=
\sum_{j,k} a_j^* q_{jk} b_k$ when $\xi = (\sum_jq_{ij}
a_j)_{i=1}^m$, $\eta = (\sum_kq_{ik} b_k)_{i=1}^m$.

In the nonunital case, this is necessarily more subtle as the
elements of $\HH_\infty=\cap_m\mbox{dom}\D^m$ must also satisfy
integrability conditions. In \cite{R1}, the notion of smooth module
was introduced for nonunital algebras which are local. As we are
dealing with quasi-local algebras, most of the results on smooth
modules in \cite{R1} are not applicable.

We take the attitude that:

Point~1) $\HH_\infty$ should be a continuous $\A$-module,
\vspace{-6pt}

Point~2) $\HH_\infty$ should embed continuously as a dense subspace
in the $C^*$-$A$-module $X_A = \overline{\HH_\infty}$,\vspace{-6pt}

Point~3) $X$ should be the completion of $qA^N$ for some $N$ and
some projection $q$ in $M_N(A_b)$ where $A_b$ is a unitization of
$A$, \vspace{-6pt}

Point~4) the Hermitian product $\HH_\infty\ni x,y\to x^*y$ should
have range in $\A$ (acting on the right).

Point 1) is implied by the condition of regularity.
\begin{proof}
For $x\in\HH_\infty$ and $a\in\A$ we have
\begin{align}\Vert\D^n(xa)\Vert_\HH=\Vert\dd^n(xa)\Vert_\HH&=
\Big\Vert\sum_{j=0}^n
\left(\begin{array}{c}n\\j\end{array}\right)\delta^{n-j}(a^{op})\dd^j(x)\Big\Vert_\HH\nno
&\leq \sum_{j=0}^n\left(\begin{array}{c}n\\j\end{array}\right)
\Vert\delta^{n-j}(a^{op})\Vert\,\Vert\dd^j(x)\Vert_\HH.\end{align}
The continuity of the action of $\A$ on $\HH_\infty$ now follows
easily.
\end{proof}

Point~2 above is included to ensure that we can recover the `module
of continuous sections vanishing at infinity' from $\HH_\infty$, and
it is a nontrivial condition as we shall see. Once we have a
continuous embedding, the image will be dense for our graph algebra
examples, since $A_c\subset\HH_\infty$.

Once we can recover the module $X$, we demand that it be `finitely
generated and projective' in the sense of 3): see also \cite[Theorem
8]{R1}. The examples arising from graph algebras have $A$ dense in
$X$, so taking $N=1$ and $q=id_{A_b}$ in any unitization $A_b$ of
$A$ shows that  3) is always satisfied for the gauge spectral triple
of a graph algebra.

All four points are satisfied in the unital case, so we will ignore
the case of a single loop in the following, focussing attention
instead on the directed trees.

Roughly speaking, without points 2)~and~4), $\HH$ can contain many
`functions' on the graph which are unbounded, and so are not in the
algebra $A$ or the module $X$. Modules of unbounded `functions' are
not terrible per se, but we prefer to remain close to the
$C^*$-theory.

{\bf Example} Let $E$ be the `dyadic directed tree'
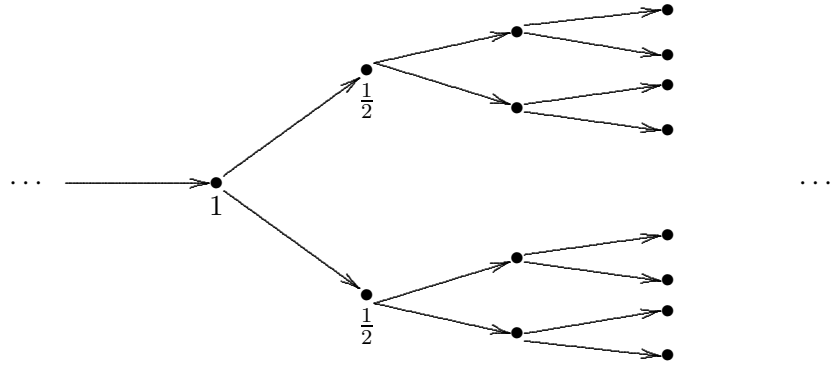
\[
\beginpicture
 \setcoordinatesystem units <1cm,1cm>
 \setplotarea x from 0 to 12, y from -3.5 to 3.5
 \put{$\cdots$} at 0.5 0
 \put{$\bullet$} at 3 0
 \put{$1$} at 3 -0.3
 \put{$\bullet$} at 5 1.5
 \put{$\frac12$} at 5 1.1
 \put{$\bullet$} at 5 -1.5
 \put{$\frac12$} at 5 -1.9
 \put{$\bullet$} at 7 2
 \put{$\bullet$} at 7 1
 \put{$\bullet$} at 7 -1
 \put{$\bullet$} at 7 -2
 \put{$\bullet$} at 9 2.3
 \put{$\bullet$} at 9 1.7
 \put{$\bullet$} at 9 1.3
 \put{$\bullet$} at 9 0.7
 \put{$\bullet$} at 9 -0.7
 \put{$\bullet$} at 9 -1.3
 \put{$\bullet$} at 9 -1.7
 \put{$\bullet$} at 9 -2.3
 \put{$\cdots$} at 11 0
 \arrow <0.25cm> [0.2,0.5] from 1 0 to 2.9 0
 \arrow <0.25cm> [0.2,0.5] from 3.1 0.1 to 4.9 1.4
 \arrow <0.25cm> [0.2,0.5] from 3.1 -0.1 to 4.9 -1.4
 \arrow <0.25cm> [0.2,0.5] from 5.1 1.6 to 6.9 2
 \arrow <0.25cm> [0.2,0.5] from 5.1 1.6 to 6.9 1.05
 \arrow <0.25cm> [0.2,0.5] from 5.1 -1.6 to 6.9 -1.05
 \arrow <0.25cm> [0.2,0.5] from 5.1 -1.6 to 6.9 -2
 \arrow <0.25cm> [0.2,0.5] from 7.1 2.1 to 8.9 2.3
 \arrow <0.25cm> [0.2,0.5] from 7.1 2 to 8.9 1.7
 \arrow <0.25cm> [0.2,0.5] from 7.1 1.05 to 8.9 1.3
 \arrow <0.25cm> [0.2,0.5] from 7.1 0.95 to 8.9 0.7
 \arrow <0.25cm> [0.2,0.5] from 7.1 -0.95 to 8.9 -0.7
 \arrow <0.25cm> [0.2,0.5] from 7.1 -1.05 to 8.9 -1.3
 \arrow <0.25cm> [0.2,0.5] from 7.1 -2 to 8.9 -1.7
 \arrow <0.25cm> [0.2,0.5] from 7.1 -2.1 to 8.9 -2.3

\endpicture
\]
Define a faithful trace as follows. If $v$ is a vertex before the
first split, let $\tau(p_v)=1$. If $v$ occurs after $n$ splits and
before $n+1$ splits, define $\tau(p_v)=2^{-n}$. Finally define
$\tau(S_\mu S_\nu^*)=\delta_{\mu,\nu}\tau(p_{r(\mu)})$. Then the
Hilbert space $\HH=L^2(X,\tau)$ contains \begin{equation}
a=\lim_{N\to\infty}\sum_{i=1}^N 2^{i/4}
p_i\label{biga}\end{equation} where $\tau(p_i)=2^{-i}$, and the
$p_i$ are mutually orthogonal. The element $a \in \HH$ in equation
(\ref{biga}) does not lie in the $C^*$-module $X$, as the limit does
not exist in the norm $\Vert\cdot\Vert_X$.

\begin{ncond}[Nonunital Finiteness]
\label{cn:finite-nonunital} The dense subspace of $\HH$ which is the
smooth domain of~$\D$,
$$
\HH_\infty := \bigcap_{m\geq 1} {\rm dom}\,\D^m
$$
has a right inner product $\A$-module structure. Moreover,
$\HH_\infty$ embeds as a dense subspace of a $C^*$-$A$-module which
is finitely generated and projective over some unitization $A_b$ of
$A$.
\end{ncond}

Having identified a working generalisation of the finiteness
condition, we identify the restrictions it places on a graph
$C^*$-algebra. So to check that New
Condition~\ref{cn:finite-nonunital} holds, we must verify points
2)~to~4).

\begin{prop}\label{finprojcore}  Suppose that the locally
finite directed graph $E$ has no sinks, no loops and satisfies the
single entry condition.
The $\A$-module $\HH_\infty$ satisfies 2)  if and only if the
$K$-theory of $A$ is finitely generated. In this case the Hilbert
space $\HH$ also satisfies point~2). If the $K$-theory of $A$ is
finitely generated then point~4) holds.
\end{prop}

{\bf Remark} Thus for the directed tree examples, the finiteness
condition is satisfied if and only if the $K$-theory of $A$ is
finitely generated.

\begin{proof}
We begin with condition 2) for our directed trees. First of all we
must have $\ker\D\cap\HH_\infty=L^2(F_c,\tau)\subset X$. Thus we
require a $C>0$ such that
$$\Vert f\Vert_X^2=\Vert f^*f\Vert_F^{1/2}\leq C\tau(f^*f)^{1/2}=C\Vert f\Vert_\HH^2,$$
for all $f\in L^2(F_c,\tau)$. In particular, we require for all
$v\in E^0$ that
$$1=\Vert p_v\Vert_F\leq C\tau(p_v)^{1/2}.$$
Hence $\tau(p_v)$ must be bounded below, which implies, by the
definition of a graph trace and the faithfulness of $\tau$, that
there exist at most finitely many ends, and so $K_0(A)$ is finitely
generated. Thus the condition is necessary.

Conversely, suppose that $K_0(A)$ is finitely generated, and let
$\mbox{rank}(K_0(A))=N<\infty$ be the number of ends. Observe that
having finitely many ends implies that any faithful graph trace is
bounded from below. Then if $f\in F_c$, Lemma \ref{describeF} allows
us to write
$$ f=\sum_{v\in E^0}\sum_{n=1}^Nc_{v,n}S_{v,n}S_{v,n}^*,\ \ \ c_{v,n}\in\C$$
where $(v,n)$ denotes a path with source $v$ and range in the $n$-th
end. We have $\Vert f\Vert_F^2=\sup|c_{v,n}|^2$.

Now suppose that $f\in \HH$, so that
$$\Vert
f\Vert_\HH^2=\tau(f^*f)=\sum_v\sum_n|c_{v,n}|^2\tau(p_n)<\infty$$
where $p_n$ is any projection in the $n$-th end. Then \bean\Vert
f\Vert_\HH^2&=&\sum_v\sum_n|c_{v,n}|^2\tau(p_n)\geq
\mbox{min}\{\tau(p_n)\}\sum_v\sum_n|c_{v,n}|^2\nno
&\geq&\mbox{min}\{\tau(p_n)\}\sup_{v,n}|c_{v,n}|^2=\mbox{min}\{\tau(p_n)\}\Vert
f\Vert_F^2=\mbox{min}\{\tau(p_n)\}\Vert f\Vert_X^2.\eean Hence $f\in
X$. Finally, suppose that $x\in\HH$, so $x=\sum_{k\in\Z}x_k$ and
$\sum_{k\in\Z}\tau(x_k^*x_k)<\infty$. As $f_k:=x_k^*x_k\in F$ and is
positive, we have \bean \Vert
x\Vert_\HH^2&=&\sum_k\tau(x^*_kx_k)=\sum_k\Vert f_k\Vert_\HH\geq
(\mbox{min}\{\tau(p_n)\})^{1/2}\sum_k\Vert f_k\Vert_X\nno
&=&(\mbox{min}\{\tau(p_n)\})^{1/2}\sum_k\Vert
x^*_kx_k\Vert_X=(\mbox{min}\{\tau(p_n)\})^{1/2}\sum_k\Vert(x_k^*x_k)^2\Vert^{1/2}_F\nno
&=&(\mbox{min}\{\tau(p_n)\})^{1/2}\sum_k\Vert
x_k^*x_k\Vert_F=(\mbox{min}\{\tau(p_n)\})^{1/2}\sup_k\Vert
x^*_kx_k\Vert_F\nno &=&(\mbox{min}\{\tau(p_n)\})^{1/2}\Vert
x\Vert_X^2.\eean

This proves that the finite generation of $K_0(A)$ is necessary and
sufficient for the second point.

For point 4), we assume that $K_0(A)$ is finitely generated. We
observe that if $x,y\in X_c=A_c\subset\HH_\infty$ we have $x^*y\in
A_c\subset\A$. In particular, $x^*y$ is in the smooth domain of the
derivation $\delta=[\dd,\cdot]$. Thus for $x,y\in X_c$ we have, by
Lemma \ref{smoothquasi},
$$\Vert\delta^m(x^*y)\Vert^2\leq\sum_{k,l\in\Z}|k-l|^{2m}\Vert
x_k^*y_l\Vert^2_A$$ where the sum over $k,l$ is finite and we have
used $\Vert\delta^m((x^*y)^{*op})\Vert=\Vert\delta^m(x^*y)\Vert$ for
$a\in A_c$ to avoid writing $op$ throughout the following
calculation. Now
$$\Vert x^*_ky_l\Vert^2_A=\Vert y^*_lx_kx^*_ky_l\Vert_A\leq \Vert
y^*_ly_l\Vert_A\Vert x_kx_k^*\Vert_A\leq
C^2\tau(y^*_ly_l)\tau(x^*_kx_k),$$ the last inequality using the
finite generation of $K_0(A)$. So we have the inequality \bea
\Vert\delta^m(x^*y)\Vert&\leq&C^2\sum_{k,l\in\Z}|k-l|^{2m}\tau(y^*_ly_l)\tau(x^*_kx_k)\nno
&\leq&C^2\sum_{k,l\in\Z}(|k|+|l|)^{2m}\tau(y^*_ly_l)\tau(x^*_kx_k)\nno
&=&C^2\sum_{k,l\in\Z}\sum_{j=0}^{2m}\bca
2m\\j\eca|k|^{2m-j}|l|^j\tau(y_l^*y_l)\tau(x^*_kx_k)\nno
&=&C^2\sum_{k,l\in\Z}\sum_{j=0}^{2m}\bca
2m\\j\eca\tau((\dd^{j/2}y_l)^*(\dd^{j/2}y_l))
\tau((\dd^{(2m-j)/2}x_k)^*(\dd^{(2m-j)/2}x_k))\nno
&=&C^2\sum_{j=0}^{2m}\bca
2m\\j\eca\Vert\dd^{j/2}y\Vert^2_\HH\Vert\dd^{(2m-j)/2}x\Vert^2_\HH.\label{long}\eea
So suppose that $\{x^i\}\subset X_c$ is a sequence converging to
$x\in\HH_\infty$ in the topology determined by the seminorms
$x\to\Vert\dd^m x\Vert_\HH$, $m\geq 0$, and similarly $y^i\to y$.


The estimate (\ref{long}) shows that
\bean\Vert\delta^m(x^*_jy_j-x_i^*y_i)\Vert_A^2&=&
\Vert\delta^m(x^*_jy_j-x^*_jy_i+x^*_jy_i-x_i^*y_i)\Vert_A^2\nno
&\leq& C^2\sum_{j=0}^{2m}\bca
2m\\j\eca\Vert\dd^{j/2}(y_j-y_i)\Vert^2_\HH\Vert\dd^{(2m-j)/2}(x_j)\Vert^2_\HH\nno
&&+C^2\sum_{j=0}^{2m}\bca
2m\\j\eca\Vert\dd^{j/2}(y_i)\Vert^2_\HH\Vert\dd^{(2m-j)/2}(x_j-x_i)\Vert^2_\HH,\eean
and this goes to zero. Hence the sequence $x_j^*y_j$ is Cauchy in
$\A$, and so for the limits  $x,y\in\HH_\infty$, the inner product
$(x|y)_A=x^*y$ is in the completion of $A_c$ for the
$\delta$-topology, and so $x^*y\in\A$.
\end{proof}
{\bf Remark} We also note that Connes stipulates that when we
restrict the Hilbert space inner product to $\HH_\infty$ we should
have
$$\la x,y\ra=\bigintcross (x|y)_\A(1+\D^2)^{-1/2}$$
where the Hermitian product is the $\A$-valued one: $(x|y)_\A=x^*y$.
However, the trace satisfies $\tau=\tau\circ\Phi$, so
$$ \tau((x|y)_\A)=\tau(x^*y)=\tau(\Phi(x^*y))=\tau((x|y)_R),$$
and the inner product does indeed satisfy this formula, up to a
factor of 2; see Equation \eqref{eq:2}. The factor of 2 also occurs
in the type I case, and is simply a matter of normalisation of the
inner product, and does not affect the Hilbert space; see
\cite[Section 5]{R2} for the constants in the commutative case.

Next we have the first order condition which specifies the bimodule
structure. In the original type I setting we have

\begin{ocond}[First Order]
\label{cn:first-ord} There are commuting representations $\pi: \A
\to \B(\HH)$ and $\pi^{op} : \A^{op} \to \B(\HH)$ of the opposite
algebra $\A^{op}$ (or equivalently, an antirepresentation of~$\A$).
Writing $a$ for $\pi(a)$, and $b^{op}$ for $\pi^{op}(b)$, we ask
that $[a, b^{op}] = 0$. In addition, the bounded operators in
$[\D,\A]$ commute with $\A^{op}$; in other words,
\begin{equation}
 [[\D, a], b^{op}] = 0 \  \ {\rm for\  all}\ \   a,b \in \A.
\label{eq:first-ord}
\end{equation}
\end{ocond}
 In the
type I setting the first order condition gives us a spectral triple
for $\A\otimes\A^{op}$, but we believe this is not essential, and
just an artefact of the type I setting. Rather we focus on the fact
that in the type I setting the algebra $\cda$ is contained in the
endomorphism algebra of the right $\A$ module $\HH_\infty$.

The finiteness condition  asks that $\HH_\infty=\bigcap_{m\geq
1}\mbox{dom}\ \D^m$ be a finite projective (right) $\A$ module. The
first order condition  then says that $\cda \subseteq
End^R_\A(\HH_\infty)$, where $R$ is for right. One would expect this
finite projective condition to be symmetric in some sense, but this
is an extra requirement. If $\HH_\infty$ is also a finite projective
left $\A$-module, then ${\mathcal C}_\D(\A^{op})\subseteq
End^L_\A(\HH_\infty)$, $L$ for left. Typically however, these two
algebras of endomorphisms, one left and one right, will {\em not}
commute with each other. They do for the gauge spectral triple of a
graph algebra, but this is a one-dimensional phenomenon (see also
\cite{GGISV}).

A moment's thought shows that regarding the (sections of the) spinor
bundle of a spin manifold $M$ as a $C^\infty(M)$ bimodule, the two
collections of endomorphisms we obtain do not commute, since both
algebras of endomorphisms are the same Clifford algebra.

These arguments, together with the proof of the reconstruction
theorem in \cite{RV}, show that the most important aspect of the
first order condition is that the algebra $\cda$ acts as
endomorphisms of a noncommutative bundle, and that the `symbol' of
$\D$ is such an endomorphism.

Moreover, in the semifinite setting we begin with a representation
$\pi:\A\to\cn\subset\B(\HH)$. The von Neumann algebra $\cn$ is thus
required to contain $\A$ and the spectral projections of $\D$, and
these are the only requirements. So typically,
$\pi^{op}(\A^{op})\not\subset\cn$, and this is the case for the
gauge spectral triple. In particular, $\A^{op}$ need not lie in the
domain of the trace we employ, and even supposing we have a version
of the first order condition, we will not obtain a spectral triple
for $\A\otimes\A^{op}$.

We therefore change the first-order condition only very slightly as
follows:
\begin{ncond}[Semifinite First Order]
\label{cn:first-ord-sf} There are commuting representations $\pi: \A
\to \cn$ and $\pi^{op} : \A^{op} \to \B(\HH)$ of the opposite
algebra $\A^{op}$.  Writing $a$ for $\pi(a)$, and $b^{op}$ for
$\pi^{op}(b)$, we ask that $[a, b^{op}] = 0$. In addition, the
bounded operators in $[\D,\A]$ commute with $\A^{op}$; in other
words,
\begin{equation}
 [[\D, a], b^{op}] = 0 \  \ {\rm for\  all}\ \   a,b \in \A.
\label{eq:first-ord-sf}
\end{equation}
\end{ncond}

For the gauge spectral triple of a directed graph, the Hilbert space
naturally carries commuting representations of $\A$ and $\A^{op}$.
The first order condition
$$ [[\D,\A],\A^{op}]=0$$
follows since $[\D,\A]\subset\A$, and the left and right actions of
$\A$ on the Hilbert space commute.

The following condition describes a spin$^c$ structure for the
noncommutative manifold, \cite{P}.

\begin{ocond}[Spin$^c$]
\label{cn:pmorita} The $C^*$-$A$-module completion of $\HH_\infty$
is a Morita equivalence bimodule between $A$ and the norm completion
of the algebra $\cda$ generated by $\A$ and $[\D,\A]$.
\end{ocond}

Since for graph algebras the $A$-bimodule $A$ is contained in $X$,
we have a natural Morita equivalence bimodule between $A$ and $A$.
As the norm closed algebra generated by $\A$ and $[\D,\A]$ is just
$A$ in the case of the gauge spectral triple, the Morita equivalence
follows. Thus there is no need to alter the spin$^c$ condition to
deal with semifiniteness or lack of a unit (at least for graph
algebras).

\begin{ncond}[Spin$^c$]
The $C^*$-$A$-module completion of $\HH_\infty$ is a Morita
equivalence bimodule between $A$ and the norm completion of the
algebra $\cda$ generated by $\A$ and $[\D,\A]$.
\end{ncond}

In the case where $\A=C^\infty(M)$, $M$ a manifold, the spin$^c$
condition (together with orientability) provides a spin$^c$
structure for $M$, \cite{P}. Given a spin$^c$  manifold $M$, $M$ is
spin if and only if at least one (oriented) Morita equivalence
bimodule admits a bijective antilinear map satisfying the
requirements of the reality condition, \cite[Theorem 9.6]{GVF}. Thus
the reality condition below, in conjunction with the spin$^c$
condition, may be regarded as a noncommutative spin structure.

\begin{ocond}[Reality]
\label{cn:real} There is an antiunitary operator $J : \HH \to \HH$
such that $J a^* J^{-1} = a^{op}$ for all $a \in \A$; and moreover,
$J^2 = \pm 1$, $J \D J^{-1} = \pm\D$ and also $J \Gamma J^{-1} =
\pm\Gamma$ in the even case, according to the following table of
signs depending only on $p \bmod 8$:
\begin{equation}
\begin{array}[t]{|c|cccc|}
\hline
p \bmod 8            & 0 & 2 & 4 & 6 \rule[-5pt]{0pt}{17pt} \\
\hline
J^2 = \pm 1          & + & - & - & + \rule[-5pt]{0pt}{17pt} \\
J\D J^{-1} = \pm\D   & + & + & + & + \rule[-5pt]{0pt}{17pt} \\
J\Gamma J^{-1} = \pm\Gamma & + & - & + & - \rule[-5pt]{0pt}{17pt} \\
\hline
\end{array}
\qquad\qquad
\begin{array}[t]{|c|cccc|}
\hline
p \bmod 8          & 1 & 3 & 5 & 7 \rule[-5pt]{0pt}{17pt} \\
\hline
J^2 = \pm 1        & + & - & - & + \rule[-5pt]{0pt}{17pt} \\
J\D J^{-1} = \pm\D & - & + & - & + \rule[-5pt]{0pt}{17pt} \\
\hline
\end{array}
\belowdisplayskip=1pc \label{eq:j-signs}
\end{equation}
For the origin of this sign table in $KR$-homology, we refer
to~\cite[\S 9.5]{GVF}.
\end{ocond}

For the gauge spectral triple, the operator $J:L^2(X,\tau)\to
L^2(X,\tau)$, $J(x)=x^*$ satisfies the reality condition for $p=1$,
namely, $J^2=1$, $Ja^*J=a^{op}$ and $J\D J=-\D$, so the bimodule and
spectral triple are real. This can be directly verified with ease.

For this reason we retain the reality condition in its original
form.

\begin{ncond}[Reality]
There is an antiunitary operator $J : \HH \to \HH$ such that $J a^*
J^{-1} = a^{op}$ for all $a \in \A$; and moreover, $J^2 = \pm 1$, $J
\D J^{-1} = \pm\D$ and also $J \Gamma J^{-1} = \pm\Gamma$ in the
even case, according to the following table~\eqref{eq:j-signs} of
signs.
\end{ncond}

We now return to the conditions. For the type I case connectedness
of the underlying noncommutative space is formulated in the
following condition. This condition has thus far only been discussed
in the commutative framework.

\begin{ocond}[Irreducibility]
\label{cn:irred} The spectral triple $(\A,\HH,\D)$ is
\textit{irreducible}: that is, the only operators in $\B(\HH)$
commuting with $\D$ and all $a \in \A$ are the scalars.
\end{ocond}

In a von Neumann algebra context it is clear what we should replace
this condition with.

\begin{ncond}[Semifinite Irreducibility]
\label{cn:irred-sf} The semifinite spectral triple $(\A,\HH,\D)$ is
\textit{irreducible}: that is, the only operators in $\cn$ commuting
with $\D$ and all $a \in \A$ are the scalars.
\end{ncond}

For our algebra $\A$, only the fixed-point subalgebra $F$ commutes
with $\D$. For graphs satisfying the single entry condition, $F$ is
abelian. A graph-theoretic argument shows that if $E$ is connected,
then no nontrivial element of $F$ can commute with all of $\A$.

We summarise our results for graph algebras.

\begin{thm}\label{11/12ths} Let $E$ be a connected locally finite
  graph
with no sinks, admitting a faithful graph trace, satisfying the
single entry condition and having finitely generated $K$-theory.
Then  the gauge spectral triple $(\A,\HH,\D)$ of $E$ satisfies the
new (semifinite, nonunital) conditions 1~to~9. If $E$ is not a
single loop, the gauge spectral triple is both nonunital
(noncompact) and semifinite.
\end{thm}

\section{$k$-Graph Manifolds}\label{sec:kgraph}

In \cite{PRS} we adapted the construction of \cite{PRen} described
earlier to construct a Kasparov module and semifinite spectral
triple for suitable $k$-graph algebras. This was accomplished by
`pushing forward' the Dirac operator (of the simplest spin
structure) on the $k$-torus, using the canonical $\T^k$ action on a
$k$-graph algebra.

We will not go into the details of these constructions as we did for
graph algebras, noting only that they are essentially analogous to
the graph case. We also omit a general discussion of $k$-graph
algebras, as this is lengthy. We will adopt the definitions,
notations and conventions of \cite{PRS}, and refer the reader to
this work for an introduction to $k$-graph algebras adapted to this
context.

We do require several notational reminders so that we can state our
results here with the minimum of ambiguity. In particular:

{\bf Warning} In this section we reverse our conventions regarding
range and source of edges. This means that sinks and sources play
opposite roles, the single entry condition becomes the single exit
condition, and so on. This is in keeping with the notation employed
in \cite{PRS}.

Briefly, a $k$-graph is a set $\Lambda$ of paths with a degree map
$d : \Lambda \to \NN^k$. For $n \in \NN^k$, we write $\Lambda^n$ for
$d^{-1}(n)$, and regard $\Lambda^0$ as the set of vertices. Paths
have the unique factorisation property: if $d(\lambda) = m + n$ then
there are unique paths $\mu \in \Lambda^m$ and $\nu \in \Lambda^n$
such that $\lambda = \mu\nu$. In particular, if $m \le n \le l =
d(\lambda)$, then there is a unique factorisation $\lambda =
\lambda(0,m)\lambda(m,n)\lambda(n,l)$ where $d(\lambda(0,m)) = m$
and so forth. It also follows that each path $\lambda$ has a unique
range $r(\lambda) \in \Lambda^0$ such that $r(\lambda)\lambda =
\lambda$; likewise for sources. With this in mind, we write
$v\Lambda^n$ for $r^{-1}(v) \cap \Lambda^n$ and $\Lambda^n v$ for
$s^{-1}(v) \cap \Lambda^n$ for $n \in \NN^k$ and $v \in \Lambda^0$.

The $C^*$-algebra $C^*(\Lambda)$ of a $k$-graph $\Lambda$ is the
universal $C^*$-algebra generated by a set $\{S_\lambda : \lambda
\in \Lambda\}$ of  partial isometries satisfying Cuntz-Krieger type
relations \cite{KP}.

For the remainder of this section, `$k$-graph' shall be an
abbreviation for `locally convex, locally finite $k$-graph without
sinks, which possesses a faithful $k$-graph trace'. All the
conditions below refer to the general semifinite nonunital versions
discussed for graph algebras (with appropriate changes to the
dimensions involved where necessary).

The gauge spectral triples $(\A,\HH,\D)$ for $k$-graph algebras
satisfy the new dimension, regularity (smoothness) and absolute
continuity conditions, with dimension $k$. All this is proved in
\cite{PRS}.

The new first order condition is satisfied just as in the graph
case, and the new irreducibility condition is also satisfied if the
$k$-graph is connected.

The remaining conditions which we need to verify are the new
finiteness, orientability, closedness, Morita equivalence (spin$^c$)
and reality conditions.

In order to do this, we will need to assume that our $k$-graphs are
row-finite with no sources ($0 < |v\Lambda^n| < \infty$ for $v \in
\Lambda^0$ and $n \in \NN^k$), and satisfy the \emph{single exit
condition} ($|\Lambda^n v| = 1$ for each $v \in \Lambda^0$ and $n
\in \NN^k$).

{\bf Finiteness and Morita Equivalence}

The proof of our rather strict finiteness condition for $k$-graphs
is almost identical to the proof for the $1$-graph case. In fact,
once we have the following result it is virtually identical.

Suppose that $\Lambda$ is a row-finite $k$-graph with no sources
satisfying the single-exit condition. We claim there is an
isomorphism of the fixed point algebra $C^*(\Lambda)^\gamma$ onto
$C_0(\Lambda^\infty)$ (where the infinite path space
$\Lambda^\infty$ is endowed with the topology generated by the
cylinder sets $\lambda\Lambda^\infty$, $\lambda \in \Lambda$). The
isomorphism  takes $S_\lambda S^*_\lambda$ to the characteristic
function $\chi_{\lambda\Lambda^\infty}$ for each $\lambda \in
\Lambda$. To see this, first recall from \cite{FPS} (see also
\cite{KP}) that for an arbitrary row-finite $k$-graph $\Lambda$ with
no sources, there is an isomorphism of the diagonal $D(\Lambda)
:=\clsp\{S_\lambda S^*_\lambda : \lambda \in \Lambda\}$ onto
$C_0(\Lambda^\infty)$ which takes $S_\lambda S^*_\lambda$ to
$\chi_{\lambda\Lambda^\infty}$. We also know that
$C^*(\Lambda)^\gamma = \clsp\{S_\mu S^*_\nu : d(\mu) = d(\nu)\}$,
but the single-exit condition ensures that whenever $S_\mu S^*_\nu
\not = 0$ and $d(\mu) = d(\nu)$, we have $\mu = \nu$. Hence
$C^*(\Lambda)^\gamma = D(\Lambda)$ when $\Lambda$ satisfies the
single-exit condition, and this establishes the claim.

In particular, it is not hard to deduce from this an exact analogue
of Lemma~\ref{describeF}: if $\Lambda$ is row-finite, satisfies the
single exit condition, and has finitely many (say $N$) ends, then
each element $a$ of $F_c$ can be expressed as
\begin{equation}
a := \sum_{v \in \Lambda^0} \sum_{i=1}^N b_{(v,i)} S_{(v,i)}
S^*_{(v,i)} \label{approx}
\end{equation}
where $(v,i)$ is a path from the $i^{\rm th}$ end to $v$.

As in the graph case, there is almost nothing to prove when the
algebra is unital. This follows since then the trace of the identity
is finite, and we can compare the Hilbert space and $C^*$-module
norms easily.

For the nonunital case we have the following.

\begin{prop}\label{k-finprojcore}  Suppose that the locally
finite, locally convex $k$-graph $(\Lambda,d)$ has no sources and
satisfies the single exit condition.
The $\A$-module $\HH_\infty$ embeds continuously in the
$C^*$-$A$-module completion if and only if the $K$-theory of $A$ is
finitely generated, in this case the Hilbert space $\HH$ does also.
If the $K$-theory of $A$ is finitely generated then the $C^*$-inner
product restricted to $\HH_\infty$ takes values in $\A$.
\end{prop}

Apart from the above result describing the fixed-point algebra of
$C^*(\Lambda)$, we also require the $K$-theory computations of
\cite{PRS} which show that in the situation of
Proposition~\ref{k-finprojcore} the $K$-theory is finitely generated
if and only if there are finitely many ends. With these results in
hand, our corresponding proof for $1$-graphs can be applied with
minor modifications.

The Morita equivalence condition is now simple, since $\cda\cong \C
liff(\R^k)\otimes\A$ (or $\C liff^+(\R^k)$ in odd dimensions) and
$\HH_\infty=\A^{2^{[k/2]}}$. So the gauge spectral triple of a
$k$-graph is spin$^c$.

{\bf Orientation}

In what follows we write $\One_k$ for $(1, \dots, 1) \in \NN^k$. We
denote the group of permutations of $\{1, \dots, k\}$ by $\Sigma_k$.

Fix a $k$-graph $\Lambda$ and a path $\mu \in \Lambda^{\One_k}$.
Given a permutation $\sigma \in \Sigma_k$, the factorisation
property guarantees that there is a unique factorisation
\[
\mu = \mu^\sigma_1 \mu^\sigma_2 \dots \mu^\sigma_k
\]
such that $\mu^\sigma_i \in \Lambda^{e_{\sigma(i)}}$ for $1 \le i
\le k$.

For example, let $k = 2$ and let $\mu = ef = ab$ be a commuting
square, so that $d(e) = d(b) = e_1$ and $d(f) = d(a) = e_2$. There
are two elements of $\Sigma_2$, namely the flip $(1, 2)$ and the
identity $\Id$. We have $\mu^{(1, 2)}_1 = a$ and $\mu^{(1, 2)}_2 =
b$ whilst $\mu^{\Id}_1 = e$ and $\mu^{\Id}_2 = f$.

We use the notation $(-1)^\sigma$ for the canonical homomorphism
$\sigma \mapsto (-1)^\sigma$ from $\Sigma_k$ to $\{-1,1\}$ which
takes the $2$-cycles $(i, j)$ to $-1$.

\begin{prop}
Let $\Lambda$ be a row-finite $k$-graph with no sources, and suppose
that for every $v \in \Lambda^0$ and $1 \le i \le k$ we have
$|\Lambda^{e_i} v| = 1$ (single exit). Define
\begin{equation}\label{eq:c def}
c_k := i^{\lceil \frac{k+1}{2} \rceil} \sum_{\mu \in
\Lambda^{\One_k}} \frac{1}{k!} \sum_{\sigma\in \Sigma_k} (-1)^\sigma
S^*_\mu \otimes S_{\mu^\sigma_1} \otimes S_{\mu^\sigma_2} \otimes
\cdots \otimes S_{\mu^\sigma_k}.
\end{equation}
Then $b(c_k) = 0$, where $b$ is the Hochschild boundary operator,
and $\pi_{\D}(c_k) = \Gamma$ where $\Gamma$ is the grading for $k$
even, and the identity for $k$ odd.
\end{prop}
\begin{proof}
We begin by establishing that $\pi_\D(c_k) = \Gamma$ because this is
the easier of the two calculations. To see this, we just calculate
(here, the $\gamma^j$ are the generators of $\C liff(\R^k)$):
\begin{align*}
\pi_D(c) &= i^{\lceil \frac{k+1}{2} \rceil} \sum_{\mu \in
\Lambda^{\One_k}} \frac{1}{k!} \sum_{\sigma\in \Sigma_k} (-1)^\sigma
S^*_\mu [\D,S_{\mu^\sigma_1}][\D,S_{\mu^\sigma_2}] \dots
[\D,S_{\mu^\sigma_k}] \\
&= i^{\lceil \frac{k+1}{2} \rceil} \sum_{\mu \in \Lambda^{\One_k}}
\frac{1}{k!} \sum_{\sigma\in \Sigma_k} (-1)^\sigma S^*_\mu
S_{\mu^\sigma_1} \gamma^{\sigma(1)} S_{\mu^\sigma_2}
\gamma^{\sigma(2)} \dots S_{\mu^\sigma_k} \gamma^{\sigma(k)} \\
&= i^{\lceil \frac{k+1}{2} \rceil} \gamma^1\cdots\gamma^k\sum_{\mu
\in \Lambda^{\One_k}} \frac{1}{k!} \sum_{\sigma\in \Sigma_k}
S^*_\mu S_{\mu^\sigma_1}  S_{\mu^\sigma_2}
\dots S_{\mu^\sigma_k}  \\
&=  \omega_{\C} \sum_{\mu \in \Lambda^{\One_k}} p_{s(\mu)},
\end{align*}
where $\omega_\C$ is the complex volume form in $\C liff(\R^k)$. The
single exit assumption ensures that the sum of vertex projections in
the last line has exactly one term for each vertex of $\Lambda$, and
hence converges to the identity in the multiplier algebra of
$C^*(\Lambda)$, establishing that $\pi_\D(c_k) = \Gamma$.

Now we need to establish that $b(c_k) = 0$. To begin with, fix $\mu
\in \Lambda^{\One_k}$. We claim that
\begin{equation}\label{eq:first step}
\begin{split}
b\Big(\sum_{\sigma\in \Sigma_k} (-1)^\sigma S^*_\mu \otimes
S_{\mu^\sigma_1}& \otimes S_{\mu^\sigma_2} \otimes \cdots \otimes
S_{\mu^\sigma_k}\Big)\\
&= \sum_{\sigma \in \Sigma_k} (-1)^\sigma \Big(S^*_{\mu^\sigma_2
\dots \mu^\sigma_k} \otimes S_{\mu^\sigma_2}
\otimes \dots \otimes S_{\mu^\sigma_k} \\
&\qquad + (-1)^k S_{\mu^\sigma_k} S^*_{\mu^\sigma_k}
S^*_{\mu^\sigma_1 \dots \mu^\sigma_{k-1}} \otimes S_{\mu^\sigma_1}
\otimes \dots \otimes S_{\mu^\sigma_{k-1}}\Big).
\end{split}
\end{equation}

To see this, we apply the definition of the Hochschild boundary $b$
to obtain
\[
\begin{split}
b\Big(\sum_{\sigma\in S_k} (-1)^\sigma S^*_\mu \otimes
S_{\mu^\sigma_1}& \otimes S_{\mu^\sigma_2} \otimes \cdots \otimes
S_{\mu^\sigma_k}\Big)\\
&= \sum_{\sigma \in \Sigma_k} (-1)^\sigma \Big(
 S^*_{\mu^\sigma_2 \dots \mu^\sigma_k} \otimes S_{\mu^\sigma_2}
   \otimes \dots \otimes S_{\mu^\sigma_k} \\
 &\qquad + \sum^{k-1}_{j=1} (-1)^j S^*_{\mu} \otimes S_{\mu^\sigma_1}
   \otimes \dots \otimes
   S_{\mu^\sigma_j} S_{\mu^\sigma_{j+1}}
   \otimes \dots \otimes
   S_{\mu^\sigma_k} \\
 &\qquad + (-1)^k S_{\mu^\sigma_k} S^*_{\mu^\sigma_k}
   S^*_{\mu^\sigma_1 \dots \mu^\sigma_{k-1}} \otimes S_{\mu^\sigma_1}
   \otimes \dots \otimes
   S_{\mu^\sigma_{k-1}}\Big).
\end{split}
\]
To establish~\eqref{eq:first step}, it therefore suffices to show
that for $1 \le j \le k-1$, we have
\[
\sum_{\sigma \in \Sigma_k} (-1)^\sigma S^*_{\mu} \otimes
S_{\mu^\sigma_1}
   \otimes \dots \otimes
   S_{\mu^\sigma_j} S_{\mu^\sigma_{j+1}}
   \otimes \dots \otimes
   S_{\mu^\sigma_k} = 0.
\]
To see this, we fix $1 \le j \le k-1$, and note that we may
partition $\Sigma_k$ as $\Sigma_k = A_j \sqcup B_j$ where $A_j :=
\{\sigma \in \Sigma_k : \sigma(j) < \sigma(j+1)\}$ and $B_j :=
\{\sigma \in \Sigma_k : \sigma(j)
> \sigma(j+1)\}$. Let $t_j \in \Sigma_k$ be the transposition $(j,\,j+1)$. Then
$\sigma \mapsto \sigma \circ t_j$ is a bijection from $A_j$ to
$B_j$.

Hence
\begin{align*}
\sum_{\sigma \in \Sigma_k} (-1)^\sigma S^*_{\mu} & \otimes
S_{\mu^\sigma_1}
   \otimes \dots \otimes
   S_{\mu^\sigma_j} S_{\mu^\sigma_{j+1}}
   \otimes \dots \otimes
   S_{\mu^\sigma_k} \\
&= \sum_{\sigma \in A_j} \Big((-1)^\sigma S^*_{\mu}
   \otimes S_{\mu^\sigma_1}
   \otimes \dots \otimes
   S_{\mu^\sigma_j} S_{\mu^\sigma_{j+1}}
   \otimes \dots \otimes
   S_{\mu^\sigma_k} \\
&\qquad + (-1)^{\sigma\circ t_j} S^*_\mu
   \otimes S_{\mu^{\sigma\circ t_j}_1}
   \otimes \dots \otimes
   S_{\mu^{\sigma\circ t_j}_j} S_{\mu^{\sigma\circ t_j}_{j+1}}
   \otimes \dots \otimes
   S_{\mu^{\sigma\circ t_j}_k}\Big).
\end{align*}

The definition of $t_j$ guarantees that $(-1)^\sigma +
(-1)^{\sigma\circ t_j} = 0$ for all $\sigma \in A_j$, and we will
therefore have established~\eqref{eq:first step} if we can show that
for fixed $1 \le j \le k-1$ and fixed $\sigma \in A_j$, we have
\begin{equation}\label{eq:elementary tensors equal}
S^*_{\mu}
   \otimes S_{\mu^\sigma_1}
   \otimes \dots \otimes
   S_{\mu^\sigma_j} S_{\mu^\sigma_{j+1}}
   \otimes \dots \otimes
   S_{\mu^\sigma_k}
= S^*_\mu
   \otimes S_{\mu^{\sigma\circ t_j}_1}
   \otimes \dots \otimes
   S_{\mu^{\sigma\circ t_j}_j} S_{\mu^{\sigma\circ t_j}_{j+1}}
   \otimes \dots \otimes
   S_{\mu^{\sigma\circ t_j}_k}
\end{equation}
By definition of $t_j$ we have $\mu^\sigma_i = \mu^{\sigma\circ
t_j}_i$ whenever $i \not= j,j+1$. If we set $m := \sum^{j-1}_{i=1}
e_{\sigma(i)} \in \NN^k$, then the factorisation property in
$\Lambda$ ensures that
\[
\mu^\sigma_j \mu^\sigma_{j+1}
 = \mu(m, m + e_{\sigma(j)} + e_{\sigma(j+1)})
 = \mu(m, m + e_{\sigma \circ t_j(j)} + e_{\sigma\circ t_j(j+1)})
 = \mu^{\sigma\circ t_j}_j \mu^{\sigma\circ t_j}_{j+1}.
\]
It follows that corresponding terms in the elementary tensors on
either side of~\eqref{eq:elementary tensors equal} are identical.
This establishes~\eqref{eq:first step}.

We must now show that if we sum the right-hand side
of~\eqref{eq:first step} over all $\mu \in \Lambda^{\One_k}$, we
obtain zero.

Fix, for the time being, $\mu \in \Lambda^{\One_k}$ and $\sigma \in
\Sigma_k$. Consider the expression
\begin{equation}\label{eq:mu term}
(-1)^\sigma S^*_{\mu^\sigma_2 \dots \mu^\sigma_k} \otimes
S_{\mu^\sigma_2} \otimes \dots \otimes S_{\mu^\sigma_k}
\end{equation}
appearing as a summand in the first term on the right-hand side
of~\eqref{eq:first step}. Let $\lambda := \mu^\sigma_2 \mu^\sigma_3
\dots \mu^\sigma_k$, so that $\mu = \mu^\sigma_1 \lambda$. Let
$\psi_k \in \Sigma_k$ be the permutation defined by $\psi_k(i) =
i+1$ for $i \le k-1$ and $\psi_k(k) = 1$. Fix $\alpha \in s(\lambda)
\Lambda^{e_{\sigma(1)}}$. Then $\lambda\alpha \in \Lambda^{\One_k}$.
Consider the expression
\[
x(\lambda,\sigma\circ\psi_k,\alpha) := (-1)^{\sigma\circ\psi_k}
(-1)^k S_{(\lambda\alpha)^{\sigma\circ\psi_k}_k}
S^*_{(\lambda\alpha)^{\sigma\circ\psi_k}_k}
S^*_{(\lambda\alpha)^{\sigma\circ\psi_k}_1 \dots
(\lambda\alpha)^{\sigma\circ\psi_k}_{k-1}} \otimes
S_{(\lambda\alpha)^{\sigma\circ\psi_k}_1} \otimes \dots \otimes
S_{(\lambda\alpha)^{\sigma\circ\psi_k}_{k-1}}
\]
which appears in the second term on the right-hand side
of~\eqref{eq:first step} for $\lambda\alpha \in \Lambda^{\One_k}$
and $\sigma\circ\psi_k \in \Sigma_k$. We have $(-1)^{\psi_k} =
(-1)^{k-1}$, and hence $(-1)^{\sigma\circ\psi_k} (-1)^k =
-(-1)^\sigma$. By definition of $\psi_k$, we have
$(\lambda\alpha)^{\sigma\circ\psi_k}_k = \alpha$, and
$(\lambda\alpha)^{\sigma\circ\psi_k}_j = \mu^\sigma_{j+1}$ for $1
\le j \le k-1$. Hence, we may rewrite
\[
x(\lambda,\sigma\circ\psi_k,\alpha) = -(-1)^\sigma S_\alpha
S^*_\alpha S^*_{\mu^\sigma_2 \dots \mu^\sigma_k} \otimes
S_{\mu^\sigma_2} \otimes \dots \otimes S_{\mu^\sigma_k}.
\]
By the Cuntz-Krieger relation, we have
\[
\sum_{\alpha \in s(\lambda) \Lambda^{e_{\sigma(1)}}} S_\alpha
S^*_\alpha = p_{s(\lambda)} = p_{s(\mu^\sigma_k)},
\]
and hence
\begin{equation}\label{eq:cancellation}
S^*_{\mu^\sigma_2 \dots \mu^\sigma_k} \otimes S_{\mu^\sigma_2}
   \otimes \dots \otimes S_{\mu^\sigma_k} -
   \sum_{\alpha \in s(\lambda) \Lambda^{e_{\sigma(1)}}}
   x(\mu^\sigma_2 \dots \mu^\sigma_k,\sigma\circ\psi_k,\alpha)
   = 0
\end{equation}

The single-exit condition and the unique factorisation property
guarantee that each
\[
S^*_{\mu^\sigma_2 \dots \mu^\sigma_k} \otimes S_{\mu^\sigma_2}
   \otimes \dots \otimes S_{\mu^\sigma_k}
\]
occurs exactly once in the first summand of the right-hand side
of~\eqref{eq:first step} as $\mu$ ranges over $\Lambda^{\One_k}$ and
$\sigma$ ranges over $\Sigma_k$. The factorisation property shows
that for fixed $\mu$ and $\sigma$, a term
$x(\lambda,\sigma',\alpha)$ is of the form
\[
x \otimes S_{\mu^\sigma_2} \otimes \dots \otimes S_{\mu^\sigma_k}
\]
for some $x \in C^*(\Lambda)$ only if $\sigma' = \sigma \circ
\psi_k$, $\lambda = \mu^\sigma_2 \dots \mu^\sigma_k$ and $\alpha \in
s(\lambda)\Lambda^{e_{\sigma(1)}}$.

Hence we may formally rewrite
\begin{align*}
b(c_k)
 &= \sum_{\mu \in \Lambda^{\One_k}} \sum_{\sigma \in \Sigma_k}
 (-1)^\sigma \Big( S^*_{\mu^\sigma_2 \dots \mu^\sigma_k}
 \otimes S_{\mu^\sigma_2} \otimes \dots \otimes S_{\mu^\sigma_k}
  - \sum_{\alpha \in s(\mu^\sigma_k)\Lambda^{e_{\sigma(1)}}}
    x(\mu^\sigma_2 \dots \mu^\sigma_k,\sigma\circ\psi_k,\alpha)\Big),
\end{align*}
which formally collapses to zero by~\eqref{eq:cancellation}.

One can check relatively easily, using the approximate identity
$\sum_{\mu \in \Lambda^{\One_k}} S_\mu S^*_\mu$ for $C^*(\Lambda)$,
that the infinite sums involved in the definition of $c_k$ and the
formal calculations in this proof make sense in the multiplier
algebra of the $k+1$-fold tensor power of $C^*(\Lambda)$.
\end{proof}

{\bf Closedness}

To show that for all $a_1,\dots,a_k\in\A$ we have
$$
\tilde\tau_\omega(\Gamma[\D,a_1]\cdots[\D,a_k](1+\D^2)^{-k/2})=0$$
it suffices to prove the result for generators of the algebra. So
let $T_{\mu_j,\nu_j}=S_{\mu_j}S_{\nu_j}^*$, $j=1,\dots,k$, be
generators. Then
$$[\D,T_{\mu_j,\nu_j}]=\gamma(id_j)=i\sum_{m=1}^k\gamma^mn_{m,j}\,T_{\mu_j,\nu_j}$$
where $d_j=(n_{1,j},\dots,n_{k,j})$ is the degree of
$T_{\mu_j,\nu_j}$. With this notation we have
\begin{align}
&\tilde\tau_\omega(\Gamma[\D,T_{\mu_1,\nu_1}]\cdots[\D,T_{\mu_k,\nu_k}](1+\D^2)^{-k/2})\nno
&=i^p\tilde\tau_\omega\left(\Gamma(\sum_{j_1}\gamma^{j_1}n_{j_1,1})\cdots(\sum_{j_k}\gamma^{j_k}n_{j_k,k})T_{\mu_1,\nu_1}\cdots
T_{\mu_k,\nu_k}(1+\D^2)^{-k/2}\right).\label{eq:first-pass}
\end{align}

Now $\Gamma=\omega_\C\otimes 1$ where $\omega_\C$ is the
(representation of) the complex volume form in the Clifford algebra.
Since  the only products of generators of the Clifford algebra with
non-zero trace are multiples of the identity, the only surviving
terms on the right hand side of Equation \eqref{eq:first-pass} when
we expand the products are those with precisely one of each
generator $\gamma^j$. Thus

\begin{align*}
&\tilde\tau_\omega(\Gamma[\D,T_{\mu_1,\nu_1}]\cdots[\D,T_{\mu_k,\nu_k}](1+\D^2)^{-k/2})\nno
&=i^p\tilde\tau_\omega\left(\Gamma\sum_{\s\in
S_k}\gamma^{\s(1)}n_{\s(1),1}\cdots\gamma^{\s(k)}n_{\s(k),k})T_{\mu_1,\nu_1}\cdots
T_{\mu_k,\nu_k}(1+\D^2)^{-k/2}\right)\nno
&=i^{p-[(p+1)/2]}\tilde\tau_\omega\left(\sum_{\s\in
S_k}(-1)^{\s}n_{\s(1),1}\cdots n_{\s(k),k})T_{\mu_1,\nu_1}\cdots
T_{\mu_k,\nu_k}(1+\D^2)^{-k/2}\right)\nno
&=i^{p-[(p+1)/2]}\det(n_{j,k})\tilde\tau_\omega\left(T_{\mu_1,\nu_1}\cdots
T_{\mu_k,\nu_k}(1+\D^2)^{-k/2}\right)
\end{align*}

Now, the trace $\tilde\tau_\omega(T_{\mu_1,\nu_1}\cdots
T_{\mu_k,\nu_k}(1+\D^2)^{-k/2})=\tau(T_{\mu_1,\nu_1}\cdots
T_{\mu_k,\nu_k})$ is zero unless $T_{\mu_1,\nu_1}\cdots
T_{\mu_k,\nu_k}\in F$, since $\tau$ is gauge invariant. This is
equivalent to
$$\sum_{j=1}^kd_j=0\Leftrightarrow\ \forall \,l\ \
\sum_{m=1}^kn_{l,m}=0.$$ Hence the first, say, column of the matrix
$(n_{j,k})$ is a linear combination of the other columns, and
$\det(n_{j,k})=0$. Hence for any generators
$T_{\mu_j,\nu_j}=S_{\mu_j}S_{\nu_j}^*$, we have
$$\tilde\tau_\omega(\Gamma[\D,T_{\mu_1,\nu_1}]\cdots[\D,T_{\mu_k,\nu_k}](1+\D^2)^{-k/2})=0.$$

{\bf Reality}

We take the complex Clifford algebra $\C liff_k$ to be generated by
$k$ elements $\gamma^j$, $j=1,\dots,k$ such that
$(\gamma^j)^*=-\gamma^j$ and
$$\gamma^j\gamma^l+\gamma^l\gamma^j=-2\delta_{l,j}Id.$$
We make some further specifications on the generators consistent
with these conventions. Denote by ${\bf j}$ the antilinear
operator on $X$ such that
$$ {\bf j}x={\bf j}\bca x_1\\\vdots\\ x_{2^{[k/2]}}\eca=\bca
x_1^*\\\vdots\\ x^*_{2^{[k/2]}}\eca.$$
 Let
$s(k)=[\frac{k}{2}](k+1)-k$ and label the generators of the
Clifford algebra so that
$$\bar{\gamma^j}={\bf j}\gamma^j{\bf j}
=\left\{\begin{array}{lr}(-1)^{s(k)}\gamma^j & j\ \mbox{odd}\\
(-1)^{s(k)+1}\gamma^j & j\ \mbox{even}\end{array}\right.$$ Observe
that $s(k)$ is even only when $k=4n$, so except for these
dimensions the odd generators have complex entries and are
invariant under transpose, while the even generators have real
entries and are antisymmetric. In dimensions $4n$ the situation is
of course reversed.

Let $\chi=\gamma^2\gamma^4\cdots\gamma^{2[k/2]}$ be the product of
the even generators (take $\chi=1$ when $k=1$). Since the entries
of $\chi$ are real for all $k$ (if $k=4n$ there are $2n$ factors
in $\chi$ and so the entries of $\chi$ are real) we have
$$ \bar\chi=\chi.$$

Using $(\gamma^j)^*=-\gamma^j$ we find
$\chi^*=(-1)^{[k/2]([k/2]+1)/2}\chi$. We then define
$$ J:=\chi\circ {\bf j}={\bf j}\circ \chi.$$

\begin{lemma} The operator $J$ satisfies $J^2=\epsilon$,
$J\D=\epsilon'\D J$ and for $k$ even $J\Gamma=\epsilon''\Gamma J$,
where $\epsilon,\,\epsilon',\,\epsilon''$ are given in the table in
Equation \eqref{eq:j-signs}.
\end{lemma}

\begin{proof} To check the sign $\epsilon$, one needs only $J^*J=1$
(which is straightforward) and
$$ J^*={\bf j}^*\circ \chi^*=(-1)^{[k/2]([k/2]+1)/2}{\bf
j}\circ\chi=(-1)^{[k/2]([k/2]+1)/2}J.$$ The sign can now be easily
checked.

The sign $\epsilon''$, in even dimensions, arises because ${\bf j}$
preserves the $\pm1$ eigenspace decomposition of $\omega_\C$, and so
commutes with $\omega_\C$, while
$\omega_\C\chi=(-1)^{k/2}\chi\omega_\C$.

For $\epsilon'$ this is more subtle. We require the straightforward
identity ${\bf j}\Phi_n{\bf j}=\Phi_{-n}$ which may be checked on
generators. Then we compute
\begin{align}
J\D\Phi_nJ^*&=J(\sum_ji\gamma^jn_j)\Phi_nJ^*=\chi(-i\sum_j{\bf
j}\gamma^j{\bf j}n_j){\bf j}\Phi_nJ^*\nno &=-i\chi\left(\sum_{j\
odd}(-1)^{s(k)}\gamma^jn_j+\sum_{j\
even}(-1)^{s(k)+1}\gamma^jn_j\right){\bf j}\Phi_nJ^*\nno
&=-i(-1)^{[(k+1)/2](k+2)}\sum_j\gamma^jn_jJ\Phi_nJ^*\nno
&=(-1)^{[(k+1)/2](k+2)}\D\Phi_{-n}.\end{align} Using the
orthogonality of the $\Phi_n$, for any $x\in\mbox{Dom}\D$ we have
$$J\D
J^*x=\sum_{n\in\Z^k}J\D\Phi_nJ^*x=
(-1)^{[(k+1)/2](k+2)}\sum_{n\in\Z^k}\D\Phi_{-n}x=(-1)^{[(k+1)/2](k+2)}\D
x.$$ The reader will check that the sign appearing here agrees
with the values of $\epsilon''$ in the table above.
\end{proof}

\begin{thm}\label{k-11/12ths} Let $(\Lambda,d)$ be a connected, locally convex, locally finite
  graph
with no sources, a faithful $k$-graph trace, satisfying the single
exit condition and having finitely generated $K$-theory. Then  the
gauge spectral triple $(\A,\HH,\D)$ of $E$ satisfies the
(semifinite, nonunital) Conditions 1 to 9.
\end{thm}


\begin{thebibliography}{99999}


\bibitem[BPRS]{BPRS} T. Bates, D. Pask, I. Raeburn and W.
Szyma\'nski, {\em The $C^*$-Algebras of Row-Finite Graphs}, New York
J. Math {\bf 6} (2000), 307--324.


\bibitem[CPS2]{CPS2} A. Carey, J. Phillips and F. Sukochev,
{\em Spectral Flow and Dixmier Traces}, Advances in Mathematics,
{\bf 173} (2003), 68--113.


\bibitem[CPRS1]{CPRS1} A. Carey, J. Phillips, A. Rennie and F. Sukochev, {\em The Hochschild Class of the Chern Character of
Semifinite Spectral Triples}, Journal of Functional Analysis, {\bf
213} (2004), 111--153.

\bibitem[CPRS2]{CPRS2} A. Carey, J. Phillips, A. Rennie and F.
Sukochev, {\em The Local Index Theorem in Semifinite von Neumann
Algebras I: Spectral Flow},  Advances in Mathematics {\bf 202}
(2006), 451--516.

\bibitem[C]{C} A. Connes,
        {\em Noncommutative Geometry},
        Academic Press, 1994.

        \bibitem[C1]{C1} A. Connes, {\em Gravity Coupled with
        Matter and the Foundation of Noncommutative Geometry},
        Commun. Math. Phys. {\bf 182} (1996), 155--176.

%

\bibitem[D]{Dix} J. Dixmier, {\em Von Neumann Algebras},
North-Holland, 1981.

\bibitem[FK]{FK} T. Fack and H. Kosaki, \emph{Generalised $s$-numbers of
$\tau$-measurable operators}, Pacific J. Math. {\bf 123} (1986),
269--300.

\bibitem[FPS]{FPS} C. Farthing, D. Pask and A. Sims, {\em Crossed
products by $\Z^l$ as higher rank graph $C^*$-algebras}, in
preparation.

\bibitem[GGISV]{GGISV} V. Gayral, J.M. Gracia-Bond\'{i}a,
B. Iochum, T. Sch\"{u}cker and J.C. Varilly, {\em Moyal Planes are
Spectral Triples}, Comm. Math. Phys. {\bf 246} (2004), 569--623.

\bibitem[GVF]{GVF} J. M. Gracia-Bond\'{i}a, J. C. Varilly and H.
Figueroa, {\em Elements of Noncommutative Geometry}, Birkhauser,
Boston, 2001.




\bibitem[KP]{KP} A. Kumjian, D. Pask, {\em Higher rank graph
$C^*$-algebras}, New York J. Math. {\bf 6} (2000) 1--20.

\bibitem[KPR]{kpr} A. Kumjian, D. Pask and I. Raeburn, {\em Cuntz-Krieger
algebras of directed graphs}, Pacific J. Math. {\bf 184} (1998),
161--174.

\bibitem[KPRR]{KPRR} A. Kumjian, D. Pask, I. Raeburn and J. Renault,
{\em Graphs, Groupoids and Cuntz-Krieger Algebras}, J. Funct. Anal.
{\bf 144} (1997), 505--541.

\bibitem[L]{L} E. C. Lance,
{\em Hilbert $C^*$-Modules}, Cambridge University Press, Cambridge,
1995.


\bibitem[M]{mal} A. Mallios,
{\em Topological Algebras, Selected Topics}, Elsevier Science
Publishers B.V., 1986.

\bibitem[P]{P} R. J. Plymen, {\em Strong Morita equivalence, spinors and symplectic spinors},
J. Operator Theory {\bf 16} (1986), 305--324.

\bibitem[PR]{PR} D. Pask and I. Raeburn, {\em On the K-Theory of
Cuntz-Krieger Algebras}, Publ. RIMS, Kyoto Univ., {\bf 32} (1996),
415--443.

\bibitem[PRen]{PRen} D. Pask and A. Rennie, {\em The Noncommutative
Geometry of Graph $C^*$-Algebras I: The Index Theorem}, J. Funct.
Anal. {\bf 233} (2006), 92--134.

\bibitem[PRS]{PRS} D. Pask, A. Rennie and A. Sims {\em The Noncommutative
    Geometry of $k$-Graph $C^*$-algebras}, $K$-theory, to appear.

\bibitem[R]{CBMSbk} I. Raeburn, \emph{Graph algebras}, CBMS Regional Conference Series in
Mathematics, Vol. 103, Amer. Math. Soc., 2005.


\bibitem[RW]{RW} I. Raeburn and D. P. Williams, {\em Morita Equivalence and
Continuous-Trace $C^*$-Algebras}, Math. Surveys \& Monographs,
vol. 60, Amer. Math. Soc., Providence, 1998.

\bibitem[RS]{RS} M. Reed and B. Simon,
{\em Volume I: Functional Analysis, Volume II: Fourier Analysis,
Self-Adjointness}, Academic Press, 1980.

\bibitem[RSz]{RSz} I. Raeburn and W. Szymanski, {\em Cuntz-Krieger
Algebras of Infinite Graphs and Matrices}, Trans. Amer. Math. Soc.
{\bf 356} (2004), 39--59.


\bibitem[R1]{R1} A. Rennie,
{\em Smoothness and Locality for Nonunital Spectral Triples},
$K$-theory, {\bf 28} (2003), 127--165.

\bibitem[R2]{R2} A. Rennie
{\em Summability for Nonunital Spectral Triples} $K$-theory, {\bf
31} (2004) pp 71-100

\bibitem[RV]{RV} A. Rennie and J. Varilly, {\em Reconstruction of Manifolds in Noncommutative Geometry},
math.OA/0610418.

\bibitem[RV2]{RV2} A. Rennie and J. Varilly, {\em Metrics for noncommutative manifolds},
in preparation.

\bibitem[RLL]{RLL} M. R\o rdam, F. Larsen and N. J. Laustsen, {\em An
Introduction to $K$-Theory and $C^*$-Algebras}, LMS Student Texts,
49, CUP, 2000.

\bibitem[S]{LBS} Larry B. Schweitzer,
{\em A Short Proof that $M_n(A)$ is local if $A$ is Local and
Fr\'{e}chet}, Int. J. math. {\bf 3} (1992), 581--589.


\bibitem[T]{T} Mark Tomforde, {\em Real Rank Zero and Tracial States of
$C^*$-Algebras Associated to Graphs}, math.OA/0204095 v2.
\end{thebibliography}
\end{document}